\documentclass{article}

\usepackage{hyperref}

\usepackage{hyperref}

\usepackage{amssymb}

\newcommand{\newc}{\newcommand}

\newc{\eqnoset}{\setcounter{equation}{0}}

\newcommand{\mref}[1]{(\ref{#1})}
\newcommand{\reflemm}[1]{Lemma~\ref{#1}}
\newcommand{\refrem}[1]{Remark~\ref{#1}}
\newcommand{\reftheo}[1]{Theorem~\ref{#1}}

\newcommand{\refcoro}[1]{Corollary~\ref{#1}}

\newcommand{\refsec}[1]{Section~\ref{#1}}

\newcommand{\beq}{\begin{equation}}
\newcommand{\eeq}{\end{equation}}
\newcommand{\beqno}[1]{\begin{equation}\label{#1}}

\newcommand{\barr}{\begin{array}}
\newcommand{\earr}{\end{array}}

\newc{\bearr}{\begin{eqnarray*}}
\newc{\eearr}{\end{eqnarray*}}

\newc{\bearrno}[1]{\begin{eqnarray}\label{#1}}
\newc{\eearrno}{\end{eqnarray}}

\newc{\non}{\nonumber}
\newc{\nol}{\nonumber\nl}

\newcommand{\bdes}{\begin{description}}
\newcommand{\edes}{\end{description}}
\newc{\benu}{\begin{enumerate}}
\newc{\eenu}{\end{enumerate}}
\newc{\btab}{\begin{tabular}}
\newc{\etab}{\end{tabular}}

\newtheorem{theorem}{Theorem}[section]
\newtheorem{defi}[theorem]{Definition}
\newtheorem{lemma}[theorem]{Lemma}
\newtheorem{rem}[theorem]{Remark}
\newtheorem{exam}[theorem]{Example}
\newtheorem{propo}[theorem]{Proposition}
\newtheorem{corol}[theorem]{Corollary}
\newtheorem{conj}[theorem]{Conjecture}

\newcommand{\btheo}[1]{\begin{theorem}\label{#1}}
\newc{\brem}[1]{\begin{rem}\label{#1}\em}
\newc{\bexam}[1]{\begin{exam}\label{#1}\em}
\newc{\bdefi}[1]{\begin{defi}\label{#1}}
\newcommand{\blemm}[1]{\begin{lemma}\label{#1}}
\newcommand{\bprop}[1]{\begin{propo}\label{#1}}
\newcommand{\bcoro}[1]{\begin{corol}\label{#1}}
\newcommand{\btheoc}[1]{\begin{conj}\label{#1}}
\newcommand{\etheo}{\end{theorem}}
\newc{\etheoc}{\end{conj}}
\newcommand{\elemm}{\end{lemma}}
\newcommand{\eprop}{\end{propo}}
\newcommand{\ecoro}{\end{corol}}
\newc{\erem}{\end{rem}}
\newc{\eexam}{\end{exam}}
\newc{\edefi}{\end{defi}}

\newc{\rmk}[1]{{\bf REMARK #1: }}
\newc{\DN}[1]{{\bf DEFINITION #1: }}

\newcommand{\bproof}{{\bf Proof:~~}}
\newc{\eproof}{{\vrule height8pt width5pt depth0pt}\vspace{3mm}}

\newc{\bfrac}[2]{\dspl{\frac{#1}{#2}}}

\newc{\nid}{\noindent}

\newcommand{\dspl}{\displaystyle}
\newc{\grad}{\nabla}
\newc{\Div}{\mbox{div}}
\newc{\pdt}[1]{\dspl{\frac{\partial{#1}}{\partial t}}}
\newc{\pdn}[1]{\dspl{\frac{\partial{#1}}{\partial \nu}}}
\newc{\pdNi}[1]{\dspl{\frac{\partial{#1}}{\partial \mathcal{N}_i}}}
\newc{\pD}[2]{\dspl{\frac{\partial{#1}}{\partial #2}}}
\newc{\dt}{\dspl{\frac{d}{dt}}}
\newc{\bdry}[1]{\mbox{$\partial #1$}}
\newc{\sgn}{\mbox{sign}}

\newc{\Hess}[1]{\frac{\partial^2 #1}{\pdh z_i \pdh z_j}}
\newc{\hess}[1]{\partial^2 #1/\pdh z_i \pdh z_j}

\newc{\ag}{\alpha}
\newc{\bg}{\beta}
\newc{\cg}{\gamma}\newc{\Cg}{\Gamma}
\newc{\dg}{\delta}\newc{\Dg}{\Delta}
\newc{\eg}{\varepsilon}
\newc{\zg}{\zeta}
\newc{\thg}{\theta}
\newc{\llg}{\lambda}\newc{\LLg}{\Lambda}
\newc{\kg}{\kappa}
\newc{\rg}{\rho}
\newc{\sg}{\sigma}\newc{\Sg}{\Sigma}
\newc{\tg}{\tau}
\newc{\fg}{\phi}\newc{\Fg}{\Phi}
\newc{\vfg}{\varphi}
\newc{\og}{\omega}\newc{\Og}{\Omega}
\newc{\pdh}{\partial}

\newc{\ccG}{{\cal G}}

\newc{\ii}[1]{\int_{#1}}
\newc{\iidx}[2]{{\dspl\int_{#1}~#2~dx}}
\newc{\bii}[1]{{\dspl \ii{#1} }}
\newc{\biii}[2]{{\dspl \iii{#1}{#2} }}
\newc{\su}[2]{\sum_{#1}^{#2}}
\newc{\bsu}[2]{{\dspl \su{#1}{#2} }}

\newc{\biiom}[1]{{\dspl\int_{\bdrom}~ #1 ~d\sg}}
\newc{\io}[1]{{\dspl\int_{\Og}~ #1 ~dx}}
\newc{\bio}[1]{{\dspl\int_{\bdrom}~ #1 ~d\sg}}
\newc{\bsir}{\bsu{i=1}{r}}
\newc{\bsim}{\bsu{i=1}{m}}

\newc{\iibr}[2]{\iidx{\bprw{#1}}{#2}}
\newc{\Intbr}[1]{\iibr{R}{#1}}
\newc{\intbr}[1]{\iibr{\rg}{#1}}
\newc{\intt}[3]{\int_{#1}^{#2}\int_\Og~#3~dxdt}

\newc{\itQ}[2]{\dspl{\int\hspace{-2.5mm}\int_{#1}~#2~dz}}
\newc{\mitQ}[2]{\dspl{\rule[1mm]{4mm}{.3mm}\hspace{-5.3mm}\int\hspace{-2.5mm}\int_{#1}~#2~dz}}
\newc{\mitQQ}[3]{\dspl{\rule[1mm]{4mm}{.3mm}\hspace{-5.3mm}\int\hspace{-2.5mm}\int_{#1}~#2~#3}}

\newc{\mitx}[2]{\dspl{\rule[1mm]{3mm}{.3mm}\hspace{-4mm}\int_{#1}~#2~dx}}
\newc{\mitmu}[2]{\dspl{\rule[1mm]{3mm}{.3mm}\hspace{-4mm}\int_{#1}~#2~d\mu}}
\newc{\iidmu}[2]{\iidx{#1}{#2}}

\newc{\iidm}[3]{{\dspl\int_{#1}~#2~d #3}}

\newc{\itQmu}[2]{\dspl{\int\hspace{-2.5mm}\int_{#1}~#2~d\mu}}
\newc{\mitQmu}[2]{\dspl{\rule[1mm]{4mm}{.3mm}\hspace{-5.3mm}\int\hspace{-2.5mm}\int_{#1}~#2~d\mu}}

\newc{\mitQq}[2]{\dspl{\rule[1mm]{4mm}{.3mm}\hspace{-5.3mm}\int\hspace{-2.5mm}\int_{#1}~#2~d\bar{z}}}
\newc{\itQq}[2]{\dspl{\int\hspace{-2.5mm}\int_{#1}~#2~d\bar{z}}}

\newc{\pder}[2]{\dspl{\frac{\partial #1}{\partial #2}}}

\newc{\bdrom}{\bdry{\Og}}

\newc{\bilhom}{\mbox{Bil}(\mbox{Hom}(\RR^{nm},\RR^{nm}))}
\newc{\VV}[1]{{V(Q_{#1})}}

\newc{\ccA}{{\mathcal A}}
\newc{\ccB}{{\mathcal B}}
\newc{\ccC}{{\mathcal C}}
\newc{\ccD}{{\mathcal D}}
\newc{\ccE}{{\mathcal E}}
\newc{\ccH}{\mathcal{H}}
\newc{\ccF}{\mathcal{F}}
\newc{\ccI}{{\mathcal I}}
\newc{\ccJ}{{\mathcal J}}
\newc{\ccK}{{\mathcal K}}
\newc{\ccP}{{\mathcal P}}
\newc{\ccQ}{{\mathcal Q}}
\newc{\ccR}{{\mathcal R}}
\newc{\ccS}{{\mathcal S}}
\newc{\ccT}{{\mathcal T}}
\newc{\ccX}{{\mathcal X}}
\newc{\ccY}{{\mathcal Y}}
\newc{\ccZ}{{\mathcal Z}}

\newc{\bb}[1]{{\mathbf #1}}

\newc{\myprod}[1]{\langle #1 \rangle}
\newc{\mypar}[1]{\left( #1 \right)}

\newc{\BLLg}{\mathbf{\LLg}}

\newc{\mA}{\mathbf{A}}
\newc{\mB}{\mathbf{B}}
\newc{\mC}{\mathbf{C}}
\newc{\mD}{\mathbf{D}}
\newc{\mE}{\mathbf{E}}
\newc{\mF}{\mathbf{F}}
\newc{\mJ}{\mathbf{J}}
\newc{\mG}{\mathbf{G}}
\newc{\mP}{\mathbf{P}}
\newc{\mR}{\mathbf{R}}
\newc{\mQ}{\mathbf{Q}}
\newc{\mX}{\mathbf{X}}
\newc{\muu}{\mathbf{u}}
\newc{\mvv}{\mathbf{v}}

\newc{\mllg}{\mathbb{\lambda}}
\newc{\mLLg}{\mathbf{\LLg}}

\newc{\lspn}[2]{\mbox{$\| #1\|_{\Lsp{#2}}$}}
\newc{\Lpn}[2]{\mbox{$\| #1\|_{#2}$}}
\newc{\Hn}[1]{\mbox{$\| #1\|_{H^1(\Og)}$}}

\newc{\mynorm}[2]{\| #1\|_{#2}}

\newcommand{\RR}{{\rm I\kern -1.6pt{\rm R}}}

\newc{\itQQ}[2]{\dspl{\int_{#1}#2\,dz}}
\newc{\mmitQQ}[2]{\dspl{\rule[1mm]{4mm}{.3mm}\hspace{-4.3mm}\int_{#1}~#2~dz}}
\newc{\MmitQQ}[2]{\dspl{\rule[1mm]{4mm}{.3mm}\hspace{-4.3mm}\int_{#1}~#2~d\mu}}

\newc{\MUmitQQ}[3]{\dspl{\rule[1mm]{4mm}{.3mm}\hspace{-4.3mm}\int_{#1}~#2~d#3}}
\newc{\MUitQQ}[3]{\dspl{\int_{#1}~#2~d#3}}

\newc{\mccP}{\mathbb{P}}
\newc{\mccK}{\mathbb{K}}

\newc{\DKTmU}{\mccK(U)}
\newc{\DKTmUold}{(K_U(U)^{-1})^T}

\newc{\myPi}{\mathbf{W}}
\newc{\myIbar}{\bar{\ccI}_1}
\newc{\myIhat}{\hat{\ccI}_1}
\newc{\myIbreve}{\breve{\ccI}_0}

\newc{\mmk}{\mathbf{k}}

\newcommand{\ma}{\mathbf{a}}

\newc{\mfu}{\mathbf{f_u}}
\newc{\mh}{\mathbf{h}}
\newc{\mb}{\mathbf{b}}
\newc{\mf}{\mathbf{f}}

\newcommand{\barrl}[2]{\barr{ll}\lefteqn{#1}\hspace{#2}&\\}

\newc{\twomatrix}[1]{\left[\barr{cc}#1\earr\right]}
\newc{\threematrix}[1]{\left[\barr{ccc}#1\earr\right]}

\newc{\mN}{\mathbf{N}}
\newc{\mI}{\mathbf{I}}
\newc{\mH}{\mathbf{H}}

\newc{\mk}{\mathbf{k}}
\newc{\mr}{\mathbf{r}}

\newc{\DIAGM}[2]{\left[\barr{ccc}#1&0\ldots&0\\
	\vdots&\ddots&\vdots\\0&\ldots0&#2\earr \right]}
\newc{\DiagM}[2]{\mbox{diag}\left[#1
	\cdots #2 \right]}
\newc{\vVEC}[2]{\left[\barr{c}#1\\
	\vdots\\#2\earr \right]}
\newc{\hVEC}[2]{\left[#1
	\cdots #2 \right]}

\newc{\mq}{\mathbf{q}}

\newc{\msys}[1]{\left\{\barr{l}#1\earr
	\right.}
\newc{\msysa}[1]{\left\{\barr{ll}#1\earr
	\right.}

\newc{\bbM}{\mathbb{M}}
\newc{\mat}[1]{\left[\barr{cc}#1\earr\right]}

\newc{\me}{\mathbf{e}}

\newc{\vecc}[2]{\left[\barr{cc}#1\\#2\earr\right]}
\newc{\mL}{\mathbb{L}}

\newc{\cO}{{\cal O}}

\newc{\cM}{{\cal M}}

\newc{\myega }{\eg_0(R)}
\newc{\myeg}{\eg_1(\eg_*)}
\newc{\myegp}{\hat{\eg}_1(\eg_*)}

\newc{\diagA}{\mathbb{A}_d}
\newc{\mBB}{\mathbb{B}}
\newc{\MLT}[1]{{\cal M}_{lt}(\Og,#1)}
\newc{\ALT}[1]{{\cal A}_{l}(\Og,#1)}
\newc{\mM}{\mathbb{M}}

\newc{\diag}[1]{\mbox{diag}(#1)}
\newc{\off}[1]{\mbox{offdiag}(#1)}
\newc{\mT}{\mathbb{T}}

\usepackage{amssymb}

\begin{document}

\vspace*{-.8in}
\begin{center} {\LARGE\em Some Maximum Principles for Cross Diffusion Systems.}

 \end{center}

\vspace{.1in}

\begin{center}

{\sc Dung Le}{\footnote {Department of Mathematics, University of
Texas at San
Antonio, One UTSA Circle, San Antonio, TX 78249. {\tt Email: Dung.Le@utsa.edu}
}}

\end{center}

\begin{abstract}
We establish  certain maximum principles for a class of strongly coupled elliptic (or cross diffusion) systems of $m\ge2$ equations. The reaction parts can be non cooperative. These new results will be crucial in obtaining coexistence and persistence for many models with cross diffusion effects.\end{abstract}

\section{Introduction}\label{intoMAX}\eqnoset

Maximum principles are important tools in analysing questions in partial differential systems such as: existence, uniqueness, positivity and symmetry of
solutions, strong positivity of associated operators, among other things. The reader may consult the classic book by Protter and
Weinberger \cite{PW} which contains most of the relevant results on this subject for partial differential equations up to
the mid-sixties.

More recently de Figueiredo and Mitidieri \cite{FM}, Sweers \cite{Sw},  L\'opez-G\'omez and  Molina-Meyer \cite{lopez}, greatly extended the theory 
to elliptic partial differential systems on a smooth bounded domain $\Og\subset \RR^N$ for $m$ equations ($N,m\ge 2$)

\beqno{sysintro}\left\{\barr{ll}-\Div(\mA DW)+\mB DW +kW-KW=F& \mbox{in $\Og$,}\\ W=0& \mbox{on $\partial \Og$,}\earr\right.\eeq   where $W=[u_1,\ldots,u_m]^T$, $\mA,\mB,K$ are $m\times m$ matrices, $k\in \RR$ and $F\in\RR^m$. The entries of these matrices (vectors) are sufficiently regular  functions/vectors (say, in $C^\nu(\Og,\RR^m)$ for some $\nu>0$  \cite{lopez}).

As usual, we assume the normal ellipticity (see \cite{Am2} and note that we do not assume $\mA$ to be symmetric): for some positive constants $\llg_0,\LLg_0$
$$\llg_0|\zeta|^2\le \myprod{\mA\zeta,\zeta}\le \LLg_0|\zeta|^2 \quad\forall\zeta\in \RR^m.$$

The main results in forementioned literature say that if $K$ is a {\em constant cooperative} matrix then \mref{sysintro}
enjoys a maximum principle provided $\mA,\mB$ are {\em diagonal} matrices. That is, if $F>0$ then $W\gg 0$ on $\Og$.
The diagonality of $\mA,\mB$ is a crucial assumption in these works. Here, we follow the standard notation: A function $f=[f_i]_1^m$ is said to satisfy $f>0$ (respectively, $f\gg 0$) if $f_i\ge 0$ for all $i$ and $f_i> 0$ for some (respectively, all) $i$'s.

In the last few decades, there is a great deal of  interest in the study of strongly coupled systems (e.g. see \cite{Am2,dlebook,SKT}). Coexistence and persistence results are investigated in \cite[Chapter 7]{dlebook1}, \cite{dlecoexper} and   they rely heavily on the spectral radii of operators associated to the linearizations of considered systems at steady states, which can be written as (for some matrix $\mathbb{G}$)
\beqno{sysintro1}\left\{\barr{ll}-\Div(\mA DW)+\mB DW +kW-KW=\mathbb{G} W& \mbox{in $\Og$,}\\ W=0& \mbox{on $\partial \Og$.}\earr\right.\eeq

The famous Krein-Rutman theorem plays an important role in the study of these spectral radius of the inverse operator associated to \mref{sysintro1} and, of course, one needs that the operator is strongly positive. Thus, certain maximum principles for \mref{sysintro} will be essential tools. 

Of course, the above maximum principles in \cite{FM,Sw,lopez} are not applicable for our purposes here because $\mA,\mB$ are full matrices. In this paper we will report some new maximum principle results concerning \mref{sysintro} where  $\mA,\mB$ can be {\em non diagonal} and $K$ can be (appropriately) {\em non cooperative}. Actually, by maximum principles in this paper we will prove that if $F$ belongs to some subcones of $\RR^m$ then $W\gg 0$. 

To the best of our knowledge, the treatment given here to the maximum principle
as well as the analysis of its relationships with the problem of the existence of
principal eigenvalues is new to all.

This paper is organized as follows. In \refsec{prelim}, we start by recalling the results in \cite{lopez} and present some counterexamples to show that they are no longer available if $\mA,\mB$ are full matrices although $K$ can be any
cooperative matrix. In \refsec{counterex}, further examples show that a change of variables will not be enough to allow us apply the result directly if $K$ is partly competitive (see also \cite{CM}). These examples, besides our main purpose, prompt us to extend the classical results.

We start \refsec{mainres} with a simple calculation for a triangular cross diffusion systems of two equations and transform it to a diagonal system in order to establish some maximum principles from the results in \cite{lopez}. We then extend the result to full cross diffusion systems like \mref{sysintro} which can be reduced to triangular form by a constant matrix. Interestingly, by allowing $\mA$ to be a full matrix we see that certain maximum principles can be proved if $K$ is appropriately competitive. We then continue to relax the condition that transformation matrix is constant. It turns out that usual maximum principles for {\em scalar} equations (we don't rely on \cite{lopez} here), an induction argument and certain assumptions on the coefficients of the system concerning the Green functions of some {\em scalar} equations should be sufficient for the proof. We will also see that the involment of Green functions seems to be necessary.

We make use of matrix notation in \refsec{matrix} to present the calculations and results of the previous section in a more compact and clear manner. It turns out that one can greatly generalize these results. In particular, we consider full cross diffusion system \mref{sysintro} when $\mA,\mB,K$ can be {\em simultaneously}  transformed (or row equivalent) to (lower or upper) triangular matrices by the same constant matrix $\mBB$, which can be non constant later.

Importantly, we give a precise and almost optimal characterization of the class $\ALT{n}$ of matrices where the argument in this work can be used to establish maximum principles. We also present a result concerning the strongly positiveness of the operator associated to \mref{sysintro}. Finally, we prsent another version of the technical \reftheo{GenMPMatTKR} in \refsec{newthm} which is easier to be verified and allows us consider the case when the reaction is completely competitive. This result reveals  important issues in applications: the cross diffusion and the reaction has to be compatible in certain ways such that a maximum principle is available; several counterexamples will be presented to support this.

Regarding \mref{sysintro} as a model in Mathematical biology/ecology, our argument reveals an interesting condition on $\mA$ for the results to be applicable : The cross diffusivities of the $i^{th}$ species must be (constant) multiples of the self diffusivity of the $j^{th}$ species if $i<j$ and, meanwhile, if $i>j$ then the $i^{th}$ species have more freedom locally.

We observe that, in the  argument, we just transformed the system into a new one, where our induction argument can work, and did not make a change of variables.   This can be generalized by combining two methods and again, we completely describe the structure of $\mA,\mB,K$ in \mref{sysintro}. We describe the conditions on such transformation matrix $\mBB$. We also give an example to show that the condition on transformation matrix $\mBB$ is necessary to conclude this paper.

\section{Preliminaries} \label{prelim}\eqnoset

Let us recall the following well known result.

\blemm{lopez} (L\'opez-G\'omez and  Molina-Meyer \cite[Theorem 3.1]{lopez}) Consider the following {\em diagonal} system for $u=[u_i]_1^m$ with homogeneous Dirichlet or Neumann boundary conditions
$$L_i(u_i) +ku_i-Ku=f_i,$$
where $L_i(\zeta)= -\Div(a_i(x)D\zeta)+b_i(x)D\zeta+c_i(x)\zeta$, a second order elliptic differential  operator with H\"older continuous coefficients. Assume that $c_0(x)+k$ is sufficiently large in terms of a given matrix $K=[k_{ij}]$ which is a cooperative matrix ($k_{ij}>0$ if $i\ne j$). Then a maximum principle holds. That is if $f_i>0$ for all $i$ then $u_i>0$ in $\Og$ for all $i$ (we can also take $K=0$ here). 

Moreover, if  $k_{ij}>0$ if $i\ne j$, then   the principal eigenvalue of $L(\fg) -K\fg=\llg_1\fg$, with $L=[L_i]_{i=1}^m$, is simple and has a positive  (vector valued) eigenfunction $\fg=[\fg_i]$.

\elemm

\brem{krem} Note that we can assume $c_0\ge0$ (by choosing $k$ large). The largeness of $c_0(x)+k$ also depends on the principal eigenfunctions of $L_i$'s. In fact, let $\llg_i,\psi_i$'s be the principal eigenpairs of $L_i$'s. $k$ should be sufficiently large such that (see the proof of \cite[Theorem 3.1]{lopez}) $$(\llg_i+k)\psi_i> \sum_j k_{ij}\psi_j.$$

\erem

\brem{lopezrem} The condition that $k_{ij}>0$ if $i\ne j$ is crucial for the last assertion on positive  eigenfunctions to hold (this is a consequence of the famous Krein-Rutman theorem for strongly positive operators, see \cite[Theorem 3.1]{lopez}; we can assert that $(L+kId-K)^{-1}$ is strongly positive if $k>0$ is sufficient large ). Moreover, it is easy to see that one can replace $k$ by a diagonal matrix $\mbox{diag}[k_1,\ldots,k_m]$ with $k_i>0$ large. In addition $k,K$ can be bounded functions (matrices) on $\Og$.
\erem

A similar result for nondiagonal systems is not true as we see in this simple counterexample.

\brem{KMP}  For any $\kappa>0$ let $\fg^{(\kappa)}$ be the positive principal eigenfunction to the positive eigenvalue $\llg_1$ of $$\left\{\barr{ll}-\Delta \fg(y)=\llg_1\fg(y)& y\in B_\kappa,\\\fg(y)=0&y\in\partial B_\kappa.\earr\right.$$  Define $u(x)=-\fg^{(\kappa)}(\kappa x), v(x)=\fg^{(\kappa)}(\kappa x)$ for $x\in B_1$. Then  $W=[u,v]^T$ is a nonpositive solution of the nondiagonal system $-\Div(\ma DW)=g$ on $B_1$ where
$$\ma=\left[\barr{cc}a&b\\b&d\earr\right],\; g=\kappa^2\llg_1\left[\barr{c} -a+b\\ -b+d\earr\right]\fg^{(\kappa)}.$$

Obviously, we can choose $a,b,d$ such that $\ma$ is elliptic (say $d>b>a>0$). One should also note that $\ma$ is  symmetric.

For any given large $k>0$ and cooperative matrix $K$ as in the lemma such that the inverse operator associated to the following problem exists \beqno{symsys}\left\{\barr{ll}-\Div(\ma DW)+kW-KW=f& \mbox{in $B_1$,}\\ W=0& \mbox{on $\partial B_1$,}\earr\right.\eeq   where 
$$ f:=g+kW-KW=\left[\barr{c} c\kappa^2\llg_1(-a+b)-k+k_{11}-k_{12}\\ c\kappa^2\llg_1(-b+d)+k+k_{21}-k_{22}\earr\right]\fg^{(\kappa)}$$

We see that $W$ is a nonpositive solution of the above system with $f$ is positive 
if we choose $\kappa$ sufficiently large. Thus, the inverse of the operator associated to \mref{symsys} exists but is not positive. The same assertion holds if $\ma$ is a upper triangular matrix. 
\erem

However, as we show later, if $b$ is not too large and $f$ belongs to certain cone then we can prove some maximum principles extending \reflemm{lopez} to establish positiveness of more general operators. In the next section, we will present some more examples which motivate us to extend the above maximum principle in this paper.

\section{Counterexamples and motivations for generalizations} \label{counterex}\eqnoset

We present here some examples which show that if the cooperative conditions are violated then maximum principles woud not hold. First of all, let $\fg$ be a positive eigenfunction to the principal eigenvalue of the problem
$$-\Delta\fg=\llg_*\fg \mbox{ in }\Og,\quad \fg=0 \mbox{ on }\partial\Og.$$

We consider the following system with $A$ is a diagonal matrix and $G=[g_{ij}]$
$$-\Div(ADW)=G(W)W \mbox{ in }\Og\times(0,\infty),\quad W=0 \mbox{ on }\partial\Og\times(0,\infty) \mbox{ and }W(x,0)=W_0(x) \mbox{ on }\Og.$$

Of course, we assume that $A$ satisfies the usual elliptic condition.

\subsection{Partial competitive case:} \label{partialcomp}  We consider the case $g_{ij}\le0$ for some $i\ne j$. 
In some cases, by  an appropriate change of variables, we can reduce this case to the cooperative ones considered previously in \reflemm{lopez}.

We assume that $\mG$ is a {\em block} matrix
$$\mG=\left[ \begin {array}{cc} A&B\\ \noalign{\medskip}C&D\end {array}
\right]
$$
where $A,D$ are square positive matrices of sizes $k,l$ and the off-diagonal entries of $B,C$ are nonnegative. We define \beqno{Pmatrix}P:=\left[ \begin {array}{cc} Id_k&0\\ \noalign{\medskip}0&-Id_l\end {array}
\right].\eeq

In some models in Mathematical Biology, it is natural to assume that the components participating in the process will react to each others in a tit for tat way. That is, the variables in $v\in\RR^{m}$ can be divided into two competing goups but they support each others in theirs owns. Thus, the symmetric (across the diagonal of $\mG$) entries will have the same sign. Using permutation matrices, we can always assume $\mG$ to have this form.

We make use of a change of variables $\bar{W}=PW$ then the  new system still satisfies the normal ellipticity. Furthermore, if $A$ is diagonal then the $P^{-1}AP$ is also diagonal and$$\bar{\mG}=P^{-1}\mG P=[\bar{g}_{ij}]=\left[ \begin {array}{cc} A&-B\\ \noalign{\medskip}-C&D\end {array}
\right].$$

Thus, the off-diagonal entries of $P^{-1}\mG P$  are nonnegative because those of $A,D,-B,-C$ are.

However, the completely competitive case is different. One can easily see that if $P$ is a diagonal matrix then it $P^{-1}\mG P$ cannot have its off-diagonal entries {\em all positive} for us to apply \reflemm{lopez}. The same problem occur if we combine $P$ with permutation matrices.

\subsection{Prey-predator case:} We consider the competitive case, i.e. $\mG$ is such that the upper diagonal entries nonnegative and the lower diagonal ones are nonpositive. We discuss only here the case $m=2$  for simplicity. Suppose that $$\mG=\left[ \begin {array}{cc} a&b\\ \noalign{\medskip}c&d\end {array}
\right],\; P=\mat{ p_{1,1}&p_{1,2}\\ p_{2,1}& p_{2,2}}
$$
where $c\le 0, b\ge0$. It is technical and quite involved but we can find $P$ such that  $P^{-1}AP$ is diagonal and make use of \reflemm{lopez}. We skip the tedious details.

\section{Main results in some simple cases} \label{mainres}\eqnoset 

We have seen that there are two limitations of the known results: 1) The systems must have diagonal main parts $A,B$; and 2) The reaction ($K$) must be cooperative. In this section we first relax the crucial assumptions that $A,B$ are diagonal and allow them to be triangular (or even full) of some special forms to obtain certain maximum principles and apply them to establish positiveness of the corresponding operators. We will present their generalizations in the next section. Let us start with a system of two equations.

\subsection{The case of two equations (a motivation)} \label{n2}

To elucidate the main ideas, we assume that $m=2$. The general case can be covered by induction. In particular, we consider the following system for $W=(u,v)$
\beqno{matsys}-\Div(ADW)+BDW+kW-KW=F\eeq
with homogeneous dirichlet boundary condition. Here, $k>0$, $K$ is a cooperative matrix and $$A=\left[\barr{cc}a&b\\c&d\earr\right],\; B=\left[\barr{cc}\bar{a}&\bar{b}\\\bar{c}&\bar{d}\earr\right],$$ with the entries are functions/vectors in $x,t$. The data $A,B,F$ are sufficiently regular (see \cite{lopez}) as usual.

For simplicity of representation we suppose that $B\equiv 0$ as the proof is similar if $B\ne0$ (that is, we allow either $A$ or $B$ are triangular/full or diagonal, as long as $a,d>0$, and $\bar{a},\bar{d}>0$ if $B\ne0$). Assume first that $A$ is {\em lower} triangular, we then consider the system of two equations on $\Og$
\beqno{matcsys} \left\{
\barr{rrr}-\Div(aDu)+ku-K_{11}u-K_{12}v&=&f_1,\\
-\Div(cDu+dDv)+kv-K_{21}u-K_{22}v&=&f_2,\earr
\right.\eeq where $k,K_{12},K_{21}>0$ and $f_1,f_2>0$ are $C^1$ functions on $\Og$. We also assume homogeneous Dirichlet boundary conditions for $u,v$ on $\partial\Og$.

Solving for $-\Div(cDu)=\frac{c}{a}(-\Div(aDu))$ ($-\Div(cDu)-\bar{c}Du$ if $B\ne0$) and assuming
that there is a constant $\cg_1\ge0$ such that (respectively, $\bar{c}=\cg_1\bar{a}$ if $B\ne0$)
\beqno{matc1}a,d>0,\;\frac{c}{a}=\cg_1,\eeq
we write the system as a ystem with diagonal main part
\beqno{matcsysa} \left\{
\barr{rrr}-\Div(aDu)+ku-K_{11}u-K_{12}v&=&f_1,\\
-\Div(dDv)+(k+\frac{c}{a}K_{12})v-K_{22}v-(K_{21}+\frac{c}{a}[k-K_{11}])u&=&f_2-\frac{c}{a}f_1,\earr
\right.\eeq

If $k, \cg_1=\frac{c}{a}\ge0$ and $k$ is large, we can apply \reflemm{lopez} (and \refrem{lopezrem}) to the above diagonal system and conclude that if $f_1,f_2-\frac{c}{a}f_1>0$ then $u,v>0$. Note that we can allow $K_{12}=0$ to obtain $u>0$ if $k$ is large as we can apply max principles for the equation of $u$ in \mref{matcsysa}. Also, if $\cg_1>0$ we can allow $K_{21}<0$ as $K_{21}+\frac{c}{a}[k-K_{11}]>0$ if $k$ is large and depends on the principal eigenfunctions of $-\Div(aD\cdot), -\Div(dD\cdot)$ (see \refrem{krem}). 

Thus, the system $$-\Div(ADW)+BDW+kW-KW=F,$$ with $A$ being lower triangular, is transformed into a diagonal system ($\bar{A}, \bar{B},\bar{k}$ are diagonal) {\em for $W$} with $\bar{K}=I_{ur}K$ being cooperative.
$$-\Div(\bar{A}DW)+\bar{B}DW+\bar{k}W-\bar{K}W=I_{ur}F,$$ where $I_{ur}$ is the (upward) row transformation matrix $I_{ur}=\left[\barr{cc}1&0\\-\cg_1&1\earr\right]$. 

Now, we consider the case $A$ is a full matrix with $b\ne0$. Our key assumption is that there are constants $\cg_1,\bg_1$ such that (same for $B$ if $B\ne0$)
\beqno{matc1}a,d>0,\;a=\cg_1 c,\; \bar{a}=\cg_1\bar{c} ,\; \frac{b}{d}=\bg_1,\;\bar{d}=\bg_1\bar{d},\; a-\frac{b}{d}c>0 \mbox{ and }\cg_1,\bg_1,\cg_1-\bg_1> 0.\eeq
These conditions imply that the entries of $A$ are nonnegative.

The same argument as above applies (solving for $-\Div(bDv)$) and 
\mref{matcsysa} becomes
$$ \left\{
\barr{rrr}-\Div([a-\frac{b}{d}c]Du)+(k+\frac{b}{d}K_{21})u-(K_{11}-\frac{b}{d}K_{21})u-(K_{12}+\frac{b}{d}k)v&=&f_1-\frac{b}{d}f_2,\\
-\Div(cDu+dDv)+kv-K_{21}u-K_{22}v&=&f_2,\earr
\right.$$

The above system is then reduced to an lower triangular system for $W$. The system becomes ($\hat{A}$ is lower triangular) by applying a downward row transformation matrix $I_{dr}=\left[\barr{cc}1&-\bg_1\\0&1\earr\right]$. 
$$-\Div(\hat{A}DW)+\hat{B}DW+\hat{k}W-\hat{K}W=I_{dr}F,$$

If we assume that  $a-\frac{b}{d}c>0$ and (the number $\cg_1^{-1}$ in \mref{matc1} for the triangular case) $(a-\frac{b}{d}c)/c=\cg_1-\bg_1>0$ then we can use the previous argument for the above triangular system again to reduce it to a diagonal system with $\bar{\hat{K}}$ cooperative
$$-\Div(\bar{\hat{A}}DW)+\bar{\hat{B}}DW+\bar{\hat{k}}W-\bar{\hat{K}}W=I_{ur}I_{dr}F,$$
and conclude that $u,v>0$. Here, $I_{ur}=\left[\barr{cc}1&0\\-\frac{1}{\cg_1-\bg_1}&1\earr\right]$.

If $\bg_1,\cg_1>0$ then we can allow $K_{12},K_{21}<0$ (of course, the numbers $K_{11}, K_{22}$ can be negative as we can choose $k>0$ large).

We summarize the above argument in the following
\btheo{maxmat} Assume that $A,B$ and $\mM=[m_{ij}]$ are full matrices at $W^*$ and \mref{matc1}. Let $\nu,k>0$ and $K=[K_{ij}]$. Define
$$L\fg=-\Div(AD\fg)+BD\fg+k\fg-K\fg,\; M\fg=\mM\fg,$$
$$C^\nu(\Og,\RR^2_F)=\{(f_1,f_2)\in C^\nu(\Og,\RR^2)\,:\,  f_1-\frac{b}{d}f_2>0\mbox{ and } f_2-\frac{c}{a}(f_1-\frac{b}{d}f_2)>0\},$$
$$C^{\nu,+}(\Og,\RR^2)=\{(f_1,f_2)\in C^\nu(\Og,\RR^2)\,:\, f_1,f_2>0\}.$$

If the entries of $k$ is large and $I_{ur}I_{dr}K$ is a cooperative matrix, then $L^{-1}:C^{\nu,+}(\Og,\RR^2_F)\to C^{\nu,+}(\Og,\RR^2)$ exists. Moreover, if the entries of the matrix $I_{ur}I_{dr}\mM$ is positive then $L^{-1}M:C^{\nu,+}(\Og,\RR^2)\to C^{\nu,+}(\Og,\RR^2)$.

Importantly, if \mref{matc1} holds then the assertions remain true if $m_{i,j}<0$ (appropriately) and $k$ is large.

\etheo

Note that $W^*$ is H\"olde continuous so that of entries of $A,B,F$ are in $C^\nu(\Og)$for some $\nu>0$.

Indeed, we proved $L^{-1}:C^\nu(\Og,\RR^2_F)\to C^{\nu,+}(\Og,\RR^2)$. The assertion on $L^{-1}M$ comes from the fact that $M:C^{\nu,+}(\Og,\RR^2)\to C^{\nu,+}(\Og,\RR^2_F)$ and this requires that $I_{ur}I_{dr}\mM$  is a positive matrix.

We apply the above result to the system $$L\fg=-\Div(AD\fg)+BD\fg=\mG\fg.$$

We have choices for $\mM$: either $\mM=\mG +kId-K$ or $\mM=\mG$. If $\mM=\mG +kId-K$ then we can allow $\mG$ is noncooperative if $K$ is. If $\mM=\mG$ then $I_{ur}I_{dr}\mM$  is a positive matrix if

$$g_{11}-\frac{b}{d}g_{21}>0 \mbox{ and } g_{12}-\frac{b}{d}g_{22}>0.$$

Note that $\bg_1,\cg_1$ depend only on $\ma^{22}$. Concerning $\mG$, we can have $g_{22}<0$ and $g_{12}<0$ if we choose $K_{12}<0$ (negatively large) and $k>0$ large. Similarly we can also have $g_{11},g_{21}<0$. This means $\mG$ can be completely competitive.

Compare to our counterexample, we can not have $L^{-1}:C^{\nu,+}(\Og,\RR^2)\to C^{\nu,+}(\Og,\RR^2)$ in general.

\subsection{The case $m>2$} \label{mbig}
We observe that, in the previous argument, we just transformed the system into a diagonal one and did not make a change of variables. This can be generalized to the case of $m$ equations. We only describe the results here and postpone their proofs to the next section where the use of matrix notation will make them much clearer.

Let $A=[a_{i,j}], B=[b_{i,j}]$, depend on $x$, be full matrices satisfying the following condition. 
\beqno{matrowconda}\barr{l}\mbox{For $i\ge1$ there are constants $\bg_i\ge0, \cg_i$ such that}\\
\left\{\barr{ll}a_{i,j}=\bg_ia_{i+1,j}&\mbox{if $j>i=1,\dots,m-1$,}\\
a_{i,j}=\cg_ia_{i-1,j}&\mbox{if $1\le j<i=2,\dots,m$},\\
\bg_{i-1}\cg_i\ne1&\forall i>1, \earr\right.
\\
\mbox{and  $a_{i,i},\;a_{i,i}-\bg_ia_{i+1,i}=a_{i,i}(1-\bg_i\cg_i)>0$ for $i\ge1$.}\earr\eeq

Here, we define $a_{m+1,j}=\cg_{m}a_{m,j}$ and $c_i=\frac{\bg_i\cg_i-1}{\bg_{i-1}-\cg_{i}^{-1}}=\cg_{i}\frac{\bg_i\cg_i-1}{\bg_{i-1}\cg_{i}-1}$. We then have
\btheo{matmaxprinciple} Assume \mref{matrowconda} and that

\beqno{bcycond}\bg_i>0\mbox{ if $1\le i<m$} \mbox{ and }\cg_i,\frac{\bg_i\cg_i-1}{\bg_{i-1}\cg_i-1}>0 \mbox{ if $1< i\le m$}.\eeq

Let $\llg_i,\psi_i$ be the principal eigenpairs of the operators $L_i\zeta=-\Div((1-\bg_i\cg_i) a_{i,i})D\zeta)$ and $K=[K_{i,j}]$ be any given matrix. Assume that $k>0$ is large such that for all $i\ge1$
$$(\llg_i+k(1+\bg_{i-1}c_i))\psi_i>\sum_j \hat{K}_{i,j}\psi_j+\sum_j k\bg_{i,j}\psi_j+\sum_j \cg_{i,j}\psi_j,$$
where $[\hat{K}_{i,j}]= [(1+c_i\bg_{i-1})K_{i,j}-\bg_iK_{i+1,j}-c_iK_{i-1,j}]$. 

Then for $k>0$ sufficiently large the inverse of the operator $L\zeta=-\Div(AD\zeta)+BD\zeta+k\zeta-K\zeta$ from $C^\nu(\Og,\RR^{m}_F)$ into $C^{\nu,+}(\Og,\RR^{m})$ exists and is positive.

\etheo

Of course, $C^{\nu,+}(\Og,\RR^{m})$ is the positive cone of $C^\nu(\Og,\RR^{m})$. The space $C^\nu(\Og,\RR^{m}_F)$ is defined by
$$C^\nu(\Og,\RR^{m}_F)=\{F\in C^{\nu,+}(\Og,\RR^{m})\,:\, \bar{I}F>0\}\mbox{ where }
\bar{I}=\mbox{diag}[1+\bg_{i-1}c_i]-[\bg]-[\cg].$$

{\bf A counterexample:} To see the importance of the eigenfunctions in our assumptions, we consider $a,b,c>0$ and the matrix
$$A=\left[\barr{ccc} a&\bg_1b&\bg_1\bg_2c\\\cg_2a&b&\bg_2 c\\\cg_3\cg_2a&\cg_3b&c \earr\right].$$

We will present positive numbers $a,b,c,\cg_{i},\bg_{i}$ such that this matrix satisfies the structure condition \mref{matrowconda} and \mref{bcycond} ($m=3$). However, we will see that \mref{bcycond} alone does not provide the desired maximum principle if there is no involvement of principal eigenfunctions of the Laplacian in our conditions.

Let $\fg$ be the positive principal eigen function of $-\Delta(\fg)=\llg_1\fg$ in $B_R$, $R>0$. We define $\psi=\fg(Rx)$ for $|x|<1$. We have $-\Delta(\psi)=\llg_1R^2\psi$. 

Let $W=[-\psi,\psi,\psi]^T$. We see that
$$-\Div(ADW)=\llg_1R^2\left[\barr{c} (-a+\bg_1b+\bg_1\bg_2c)\psi\\
(-\cg_2a+b+\bg_2 c)\psi\\
(-\cg_3\cg_2a+\cg_3b+c)\psi
\earr\right]$$

Thus, for any $k>0$ and $K=[K_{i,j}]$ the components of
$-\Div(ADW)+kW-KW$ are $$f_1=[\llg_1R^2(-a+\bg_1b+\bg_1\bg_2c)-k-(-K_{1,1}+K_{1,2}+K_{1,3})]\psi, $$
$$f_2=[\llg_1R^2(-\cg_2a+b+\bg_2c)-k-(-K_{2,1}+K_{2,2}+K_{2,3})]\psi, $$
$$f_3=[\llg_1R^2(-\cg_3\cg_2a+\cg_3b+c)-k-(-K_{3,1}+K_{3,2}+K_{3,3})]\psi. $$

Clearly, for any $k>0$ and $K$ we can achieve a counterexample to our result if \mref{bcycond} is verified but no involvement of $\psi$. We want $f_1,f_2,f_3>0$ on $B_1$ although $W<0$. Since we can choose $R>0$ large, this will be the case if
$$(i)\quad (-a+\bg_1b+\bg_1\bg_2c)>0, (ii)\quad(-\cg_2a+b+\bg_2c)>0, (iii)\quad(-\cg_3\cg_2a+\cg_3b+c)>0. $$

For simplicity, let $c=b=a>0$. (i)-(iii) are equivalent to
$$(i)\quad -1+\bg_1+\bg_1\bg_2>0, (ii)\quad-\cg_2+1+\bg_2>0, (iii)\quad-\cg_3\cg_2+\cg_3+1>0. $$
Or $(i')\quad \bg_1(1+\bg_2)>1, (ii'')\quad\bg_2>\cg_2-1, (iii''')\quad 1> (\cg_2-1)\cg_3$. Now, we choose $\bg_2\in (0,\frac{-1+\sqrt{5}}{2})$ and $\cg_2>0$ such that $1<\cg_2<\bg_2+1$. Then (ii') is satisfied and $0<\cg_2\bg_2<1$. Next, we choose $0<\bg_3=\cg_2-1$ and $\cg_3>0$ such that $\cg_3<\frac{1}{\cg_2-1}$ then (iii') is verified and $\cg_3\bg_3<1$. Thus, \mref{bcycond} is also verified. Of course, $\bg_1$ can be any positive number for (i') to hold.

\subsection{The case $\cg_i$'s are functions:} We will see that similar results still hold if $\cg_i$'s are functions we make use of the Green functions of certain scalar equations.

We assume that there are {\em constants} $\bg_i$'s with $\bg_{m}=0$ such that
\beqno{matrowcondfunct}a_{i,i}>\bg_ia_{i+1,i}\quad \forall i=1,\ldots,m,\quad a_{i,j}=\bg_ia_{i+1,j}\quad \mbox{if $j>i$}.\eeq

Also, assume that  $\mathbb{G}_j$ is the Green function of $-\Div((a_{j,j}-\bg_ja_{j+1,j})D\zeta)$ then we can prove by an easy induction in $i$ (see the next section) that if (as $(a_{j,j}-\bg_{j}a_{j+1,j})>0$) \beqno{RHScond}\left\{\barr{ll}\hat{f}_i:=f_i-\bg_if_{i+1}+\sum_{j<i}-\cg_{i,j}(f_j-\bg_jf_{j+1})>0,&\\
\myprod{D_x\cg_{i,j}(x),D_x\mathbb{G}_j(x,y)}\ge0,& \mbox{ for $j<i$ and a.e. $y$}\earr\right.\eeq then maximum principles (for scalar equations) imply that $u_i>0$ on $\Og$ for all $i$.

{\bf An example:} Consider the easiest case when $A$ is triangular and $B=0$ as in \mref{matcsys} $$\left\{
\barr{rrr}-\Div(aDu)+ku-K_{11}u-K_{12}v&=&f_1,\\
-\Div(cDu+dDv)+kv-K_{21}u-K_{22}v&=&f_2,\earr
\right.$$ where $k,K_{12},K_{21}>0$ and $f_1,f_2>0$ are $C^1$ functions on $\Og=B_1$. We also assume homogeneous boundary conditions for $u,v$ on $\partial\Og$ and $\frac{d}{a}=\dg>0$, a constant. Suppose that $\frac{c}{a} =\cg>0$ a function on $\Og$. As before, we can solve for $-\Div(cDu)=-\Div(\cg aDu)=-\cg\Div(aDu)-D(\cg)aDu$ and obtain
\beqno{contmat}\left\{\barr{lll}-\Div(aDu)+ku-K_{11}u-K_{12}v&=&f_1,\\
-\Div(dDv)+(k+\cg K_{12})v-(K_{21}-\cg k+\cg K_{11}) u-K_{22}v&=&f_2-\cg f_1+D(\cg)aDu.\earr
\right.\eeq

We wish to apply \reflemm{lopez} for diagonal systems then we want the right hand side terms and $u,v$ to be positive.

However, let $\fg>0$ be the principal eigenfunction of $-\Div(aD\fg)=\llg_1\fg$ on $B_R$. Set $\psi(x)=\fg(Rx)$, $x\in B_1$. Also, $u=\psi$, $v=-\psi$ and $\cg=\psi$. We see that
$$-\Div(aDu)+ku-K_{11}u-K_{12}v=(\llg_1+k-K_{11}+K_{12})\psi,$$
where $f_1=(\llg_1+k-K_{11}+K_{12})\psi>0$ if we choose $k$ large. Meanwhile, as $v=-\psi$,
$$\barrl{-\Div(dDv)+(k+\cg K_{12})v-(K_{21}-\cg k+\cg K_{11}) u-K_{22}v=D(\cg)aD\psi+}{1.5cm}&(-\dg\llg_1+(-\cg (\llg_1+k-K_{11}+K_{12})-(k+\cg K_{12})-(K_{21}-\cg k+\cg K_{11}) +K_{22})\psi=
\\&D(\cg)aD\psi+(-k-\dg\llg_1+(-\cg (\llg_1-K_{11}+K_{12})-\cg K_{12}-(K_{21}+\cg K_{11}) +K_{22})\psi.\earr
$$

As $D\cg=D\psi=RD\fg$. We have $D(\cg)aD\psi=aR^2|D\fg|^2$ which is large if $R,k$ are large. We see that the right hand side is positive. Thus, the right hand side terms of \mref{contmat} are positive but $W=[u,v]^T=[\psi,-\psi]^T$ is not. Note that if $\dg$ is a function then the above right hand side is
$$(D(\dg)+D(\cg))aD\psi+(k-\dg\llg_1-\cg (\llg_1-K_{11}+K_{12})+\cg K_{12}-(K_{21}+\cg K_{11}) -K_{22})\psi,$$ which is also positive if $R,k$ are large. Thus, \mref{RHScond} seems to be necessary.

\section{The general case - matrices}\label{matrix}\eqnoset
We now make use of matrices notation to simplify our previous presentation and greatly generalize our results.

\subsection{A simple case:}
We start with a simple case when one of our transformation matrices is a constant one. Define the upper triangular {\em constant} matrix 
\beqno{BBupper}
\mBB=\left[\barr{ccccc} 1&-\bg_1&0&\cdots&0\\
0&1&-\bg_2&0\cdots\\
0&0&\ddots&\ddots&0\\
\vdots&\vdots&\ddots&\ddots&-\bg_{m-1}\\
0&\cdots&\cdots&0&1\\
\earr\right]\mbox{ so that }\mBB A=[a_{i,j}-\bg_ia_{i+1,j}].\eeq

Since $\mBB A$ is a lower triangular matrix, we can let $\mL$ be the lower triangular matrix (whose entries can be functions)
$$\mathbb{L}=\left[\barr{ccccc} 1&0&\cdots&&0\\
\cg_{2,1}&1&0&\cdots&0\\
\cg_{3,1}&\cg_{3,2}&\ddots&\cdots&0\\
\vdots&\vdots&\ddots&\ddots&0\\
\cg_{m,1}&\cg_{m,2}&\cdots&\cg_{m,m-1}&1\\
\earr\right] \mbox{ such that } \mBB A=\mL^{-1}\diagA,$$
where $\diagA:=\mbox{diag}[a_{i,i}-\bg_ia_{i+1,i}]$, a diagonal matrix.

Multiplying $\mBB$ (a constant matrix) to the left of the system $-\Div(ADW)=F$, we have
$$-\Div(\mBB ADW)=\mBB F\Rightarrow -\Div(\mL^{-1}\diagA DW)=\mBB F.$$

Since $-\Div(\mL\diagA DW)=-\mL\Div(\diagA DW)-D\mL \diagA DW$, we have for $[\hat{f}_i]_1^{m}:=\mL(\mBB F+D\mL^{-1}\diagA DW)$
\beqno{diagsys}-\Div(\diagA DW)=\mL^{-1}(\mBB F+D\mL\diagA DW)\Rightarrow -\Div(\diagA DW)=[\hat{f}_i]_1^{m}.\eeq
Thus, if $\mL(\mBB F+D\mL^{-1}\diagA DW)\gg 0$ then $W\gg 0$. We also write $\mL^{-1}=[\hat{\cg}_{i,j}]$, a lower triangular matrix.

Let $\mG$ be the Green (vector) function of the diagonal operator $-\Div(\diagA D\zeta)$ (with the homogenous  Dirichlet boundary condition) then
$$-\Div(\diagA D\Fg)=\Psi\Rightarrow \Fg(x)=\int_\Og \mG(x,y)\Psi(y)~dy\Rightarrow D_x\Fg(x)=\int_\Og D_x\mG(x,y)\Psi(y)~dy.$$

Note that the diagonal entries of $\mL D\mL^{-1}$ are zeros so that the right hand side  of the $i^{th}$ equation has only terms involving $Du_j$'s for $j<i$ which can be written as
$$\sum_{j<i}\sum_{k}\myprod{\cg_{i,k}(x)D_x\hat{\cg}_{k,j}(a_{j,j}-\bg_ja_{j+1,j}),Du_j}=\sum_{j<i}(a_{j,j}-\bg_ja_{j+1,j})\sum_{k}\cg_{i,k}(x)\myprod{D_x\hat{\cg}_{k,j},\int_\Og D_x\mG_j(x,y)\hat{f}_j dy}.$$
Here, we use the fact that $-\Div((a_{j,j}-\bg_ia_{j+1,j})Du_i)=\hat{f}_j=:\Psi_j$.

We note that $\hat{f}_1>0$ does not involve with $DW$ so that if \beqno{bconstmat}\mL\mBB F\gg 0 \mbox{ and }\sum_{k}\myprod{\cg_{i,k}(x)D_x\hat{\cg}_{k,j},D_x\mG_j(x,y)}\ge0\mbox{ a.e in $x,y\in\Og$ and $\forall i>j\ge1$}\eeq then $\hat{f}_2>0$. We then use this information from the second equation in the third one to conclude that $\hat{f}_3>0$ and so on. By induction, we see that $\hat{f}_i>0$ for all $i$ so that $W\gg 0$.

We will prove the following result (where $\mBB$ is a constant matrix)
\btheo{GenMPMat} Assume that $\mA=[a_{i,j}],\mB=[b_{i,j}],K=[K_{i,j}]$ can be {\em simultaneously}  transformed (or row equivalent) to lower  triangular matrices $\mBB\mA,\mBB\mB,\mBB K$ by a {\em constant} matrix $\mBB$. Let $\mL=[\cg_{i,j}]$ be lower triangular matrix   such that $\mBB\mA=\mL^{-1}\mA_d$, $\mA_d=\mbox{diag}[\hat{a}_{1,1}\cdots,\hat{a}_{n,n}]$ with $\hat{a}_{i,i}>0$.  Suppose that $\off{\mL\mBB (K-kId)}$ is {\em lower triangular} and its off  diagonal entries of  are nonegative and $k>0$ is sufficiently large. Also, the entries of $\diag{\mL (kId-K)}$ are nonegative.

Denote $\mC=-\mL\mBB\mB + \mL D(\mL^{-1}) \mA_d=[\mathbf{c}_{i,j}]$ ({\em $\mC$ is lower triangular}). Let $\mG=\mbox{diag}[\mG_i]$ be the Green function of the diagonal system (with the homogenous Dirichlet boundary condition) $$-\Div(\mA_dD\zeta)-\diag{\mC}D\zeta+\diag{\mL (kId-K)}.$$ 

Suppose that

\beqno{Bmatcondthm} \myprod{\mathbf{c}_{i,j}(x),D_x\mG_j(x,y)}\ge0\mbox{ a.e $x,y\in\Og$ and $\forall i>j\ge1$}.\eeq

Then the inverse map ${\cal L}^{-1}$ of the map associated to the system 
$$\left\{\barr{ll} -\Div(\mA DW)+\mB DW +kW -KW=F &\mbox{in $\Og$,}\\ W=0&\mbox{on $\partial\Og$.}\earr
\right.$$
is a map on from $\{F\in C(\Og,\RR^n)\,:\, \mL\mBB F\gg 0 \}$ into $\{W\in C(\Og,\RR^n)\,:\, W\gg 0 \}$.

\etheo

\bproof Let us start with a simple case and let $\mA=[a_{i,j}], \mB=[b_{i,j}]$ be lower triangular matrices and the diagonal entries of $\mA$ are positive. We consider the following system
\beqno{gentrisys}\left\{\barr{ll} -\Div(\mA DW)+\mB DW=F &\mbox{in $\Og$,}\\ W=0&\mbox{on $\partial\Og$.}\earr
\right.\eeq

Of course, we can find a lower triangular matrix  $\mL=[\cg_{i,j}]$ such that $\mA=\mL^{-1}\mA_d=\mbox{diag}[a_{1,1}\cdots,a_{n,n}]$ with $a_{i,i}>0$. Let $\mL^{-1}=[\hat{\cg}_{i,j}]$. Because
$$-\Div(\mL^{-1}\mA_d DW)+\mB DW=F\Leftrightarrow -\mL^{-1}\Div(\mA_d DW)+\mB DW=F+D(\mL^{-1}) \mA_d DW$$ so that $ -\Div(\mA_d DW)=\mL F-\mL\mB DW+ \mL D(\mL^{-1}) \mA_d DW$.

As $a_{i,i}>0$ and the term involving $Du_i$ in the $i^{th}$ equation is $\cg_{i,i}(D\hat{\cg}_{i,i}a_{i,i}-b_{i,i})Du_i$, let $\mG_i$ be the Green function of $-\Div(a_{i,i}D\zeta)-\cg_{i,i}(D\hat{\cg}_{i,i}a_{i,i}-b_{i,i})D\zeta$ (with the homogenous  Dirichlet boundary condition).

Note that the $\mL (D\mL^{-1}\mA_d-\mB)$ is a lower triangular matrix so that the right hand side  of the $i^{th}$ equation has only terms involving $Du_j$'s for $j\le i$  and those involving $Du_j$'s for $j< i$ are, as before 
$$\barrl{\sum_{j<i}\sum_{k}\myprod{\cg_{i,k}(x)(D_x\hat{\cg}_{k,j}(x)a_{j,j}(x)-b_{k,j}(x)),D_xu_j(x)}=}{4cm}
&\sum_{j<i}\sum_{k}\myprod{\cg_{i,k}(x)(D_x\hat{\cg}_{k,j}(x)a_{j,j}(x)-b_{k,j}(x)),\int_\Og D_x\mG_j(x,y)\hat{f}_j dy}\earr$$
where $[\hat{f}_i]:=\mL F-\mL\mB DW+ \mL D(\mL^{-1}) \mA_d DW+\mbox{diag}[\cg_{i,i}(D\hat{\cg}_{i,i}a_{i,i}-b_{i,i})Du_i]$. 
Note also that if $\mB=0$ then the above is
$$\sum_{j<i}\sum_{k}\cg_{i,k}(x)\myprod{D_x\hat{\cg}_{k,j}a_{j,j}(x),Du_j(x)}=\sum_{j<i}a_{j,j}(x)\sum_{k}\cg_{i,k}(x)\myprod{D_x\hat{\cg}_{k,j}(x),\int_\Og D_x\mG_j(x,y)\hat{f}_j dy}$$
We see that $\hat{f}_i$ has only $Du_j$'s for $j<i$ in it so that the argument before the theorem can apply and shows that $W\gg 0$  if \beqno{Bmatcond}\mL F\gg 0 \mbox{ and } \sum_{k}\cg_{i,k}(x)\myprod{(D_x\hat{\cg}_{k,j}(x)a_{j,j}(x)-b_{k,j}(x)),D_x\mG_j(x,y)}\ge0\mbox{ a.e $x,y\in\Og$ and $\forall i>j\ge1$}.\eeq
Again if $\mB=0$ then the above condition is reduced to
$$\mL F\gg 0 \mbox{ and } \sum_{k}\cg_{i,k}(x)\myprod{D_x\hat{\cg}_{k,j}(x)a_{j,j}(x), D_x\mG_j(x,y)}\ge0\mbox{ a.e $x,y\in\Og$ and $\forall i>j\ge1$}.$$

When $\cg_i$'s are functions one may not be able to apply  the maximum principle in \reflemm{lopez} because we have to use the Green functions of  diagonal systems. Therefore if we include $KW$ in the discussion then it may have to be incorporated into $F$.

Now, let $\mA,\mB,\mL$ be lower triangular matrices such that $\mA=\mL^{-1}\mA_d$ for some diagonal matrix $\mA_d$. For some matrix $K$, we can consider
$$\left\{\barr{ll} -\Div(\mA DW)+\mB DW+kW -KW=F &\mbox{in $\Og$,}\\ W=0&\mbox{on $\partial\Og$.}\earr
\right.$$

To simplify the notations below, for any matrix $n\times n$ square matrix $M=[m_{i,j}]$ we write $\diag{M}=\mbox{diag}[m_{1,1},\ldots,m_{n,n}]$ and $\off{M}=M-\diag{M}$.

We add  $kW$, $k\ge0$, to both side of the system and denote $\mC=-\mL\mB+ \mL D(\mL^{-1}) \mA_d$. We write the system as
$$-\Div(\mA_d DW)-\diag{\mC} DW+\diag{\mL(kId-K)}W=\mL F+\off{\mL (KW-kW)}W+\off{\mC} DW.$$

Note that the left hand side of the system is an diagonal operator.
We then use the (vector) Green function $\mG=\mbox{diag}[\mG_i]$ of the diagonal system $-\Div(\mA_d DW)-\diag{\mC} DW+\diag{\mL(kId-K)}W$. If the matrices $\off{\mL (K-kId)}, \off{\mC}$ are lower triangular then the right hand side of the $i^{th}$ equation involves only $u_j$ for $j<i$.  So, if the entries of $\off{\mL (K-kId)}$ are nonnegative then, by induction and maximum principles for scalar equations (not for systems as in \reflemm{lopez}) to obtain $u_i>0$ so on. We then see that \mref{Bmatcond} implies $W\gg0$.

Finally, as in the statement of the theorem, assume that $\mA,\mB,K$ can be {\em simultaneously}  transformed (or row equivalent) to (lower or upper) triangular matrices by the same constant matrix $\mBB$. We simply multiply the constant matrix to the system and apply the above result with $\mA,\mB,K$ being replaced by $\mBB\mA,\mBB\mB,\mBB K$. The condition that the off  diagonal entries of $\off{\mL (K-kId)}$ are nonegative can be fulfilled if the signs of $\cg_{i,j}, K_{i,j}$ are appropriate such that \mref{Bmatcond} is verified. We see that \mref{Bmatcond} is determined by the Green function $\mG$ which in turn depends heavily on the (cross diffusion) matrices $\mA,\mB$. This allows $K$ to be non cooperative. That is, $K_{i,j}$ can be appropriately negative. The proof is complete.
\eproof

The condition \mref{Bmatcond} should also be compared with the counterexample in \refrem{KMP}.

If $\mBB\mA=[\hat{a}_{i,j}],\mBB\mB=[\hat{b}_{i,j}]$ and $\mL=[\cg_{i,j}], \mL^{-1}=[\hat{\cg}_{i,j}]$ then the assumption \mref{Bmatcondthm}  explicitly is
$$\sum_{k\le j}\cg_{i,k}(x)\myprod{(D_x\hat{\cg}_{k,j}\hat{a}_{j,j}-\hat{b}_{k,j}), D_x\mG_j(x,y)}\ge0\mbox{ a.e $y\in\Og$ and $\forall i>j\ge1$}.$$

Compare to \reftheo{matmaxprinciple}, we see that the results here (when we allow $\cg_{i,j}$'s to be functions) give us much more flexible structures of $\mA,\mB$ and choices of $K$. The number $k$ needs to be sufficiently large such that ${\cal L}^{-1}$ exists.

\brem{LopezKR} It is important that we have written the system as $$-\Div(\mA_d DW)-\diag{\mC} DW+\diag{\mL(kId-K)}W=\mL F+\off{\mL (KW-kW)}W+\off{\mC} DW$$ so that we can make use of the usual maximum principle for scalar equations. However, this does not imply the strongly positive property in order that we can use the Krein-Rutman theorem.

Instead we write the system as
$$-\Div(\mA_d DW)-\diag{\mC} DW-\mL(K-kId)W=\mL F+\off{\mC} DW$$ and we can apply \reflemm{lopez} (one can add and subtract $kId$ to the left hand side) if $\mL(K-kId)$ is cooperative and assume that $\mL F+\off{\mC} DW>0$. Under the assumption like \mref{Bmatcondthm} concerning the Green functions, we will need only that $\mL F>0$ to show that $W\gg 0$.

\erem

\subsection{$\bg_i$'s are not a constants?} One wishes to remove the assumption that $\bg_i$'s are constant. Even we can still row transform the system to a lower triangular one, if $\bg_i$'s are functions then $Du_j$'s are present in all equations so that an induction argument and an use of Green's functions as before can not be used to provide conditions to show that $W>0$. However, we will show that our argument can still be extended if $\mA,\mB$ satisfy certain (almost optimal) structural conditions.

We can allow $\mBB$ to be
a matrix (with function entries) and prove that

\bcoro{GenMPMatcoro} The assertion of \reftheo{GenMPMat} holds when $\mBB$ is a function matrix if, in addition, $D(\mBB)\mBB^{-1}$ is lower triangular.
\ecoro

\bproof We write the system as $-\Div(\mBB\mA DW)+(\mBB\mB+D(\mBB)\mA)DW=\mBB F$. Then, in order to use Green functions for diagonal systems and the lower order terms of $i^{th}$ equation are involved only with $u_j, Du_j$ for $j\le i$ (so that an induction argument can be applied), we assume that $\mL(\mBB\mB+D(\mBB)\mA)$ is lower triangular. Since $\mL,\mBB\mB$ are lower triangular we just need that $D(\mBB)\mA$ is.

As  $\mBB\mA$ is lower triangular and $D(\mBB)\mA=D(\mBB)\mBB^{-1}\mBB\mA$, we see our argument can continue if $D(\mBB)\mBB^{-1}$ is lower triangular. We easily see that we can define $\mC$ as before and use the same Green function $\mG$ and suppose that
\mref{Bmatcondthm} holds.
\eproof

Of course, \reftheo{GenMPMat} and \refcoro{GenMPMatcoro} hold if we replace the property lower triangular by upper triangular.

Let us try to understand the structure of a matrix $\mBB$ such that $D(\mBB)\mBB^{-1}$ is lower triangular.

Since we are going to use the adjective 'lower triangular' very often, we introduce the subset $\MLT{n}$ of $n\times n$ square matrices on $\Og$.

Suppose that $\ag$ is an square matrix, $\bg$ ($\cg$) is respectively column (row) vector and $\dg$ is a function. We write the {\em main diagonal} blocks of $\mBB$ in block form
$$\mBB_i=\mat{\ag&\bg\\\cg&\dg}$$

Dropping the subscript $i$, we will prove that 

\blemm{DBBlem} Assume that $\ag^{-1},\dg^{-1},(\ag-\dg^{-1}\bg\cg)^{-1}$ exist (Note that we just need $\ag^{-1},\dg^{-1}$  and $(\dg- \cg\ag^{-1}\bg)^{-1}$ exist). 

Then
$$D(\mBB)\mBB^{-1}\in\MLT{n}\Leftrightarrow D\ag\ag^{-1} \in \MLT{n-1} \mbox{ and }\bg=\ag k \mbox{ for some constant vector $k$}.$$
\elemm

We first recall the following fact, in block form \beqno{MATinverse}\mat{\ag&\bg\\\cg&\dg}^{-1}=\mat{(\ag-\frac{1}{\dg}\bg\cg)^{-1} & -\ag^{-1}\bg(\dg-\cg\ag^{-1}\bg)^{-1}\\
	-\frac{1}{\dg}\cg(\ag-\frac{1}{\dg}\bg\cg)^{-1} & (\dg-\cg\ag^{-1}\bg)^{-1}}.\eeq

\bproof
If $D(\mBB)\mBB^{-1}$ is lower triangular then from its last column we must have that 
$$D\ag[-\ag^{-1}\bg(\dg-\cg\ag^{-1}\bg)^{-1}]+D\bg(\dg-\cg\ag^{-1}\bg)^{-1}=0\Leftrightarrow[-D\ag\ag^{-1}\bg+D\bg](\dg-\cg\ag^{-1}\bg)^{-1}=0.$$
So that $\ag[-\ag^{-1}D\ag\ag^{-1}\bg+\ag^{-1}D\bg](\dg-\cg\ag^{-1}\bg)^{-1}=\ag D(\ag^{-1}\bg)(\dg-\cg\ag^{-1}\bg)^{-1}=0\Rightarrow D(\ag^{-1}\bg=0)$ and, therefore, $\ag^{-1}\bg=k$, a constant vector.

Also, we must have that the main block of $D\mBB\mBB^{-1}$ $$D\ag(\ag-\frac{1}{\dg}\bg\cg)^{-1}-D\bg\frac{1}{\dg}\cg(\ag-\frac{1}{\dg}\bg\cg)^{-1}=(D\ag-\frac{1}{\dg}D\bg\cg)(\ag-\frac{1}{\dg}\bg\cg)^{-1}$$ is lower triangular. Using the fact that $D\bg=D\ag k$, the above matrix is
$$(D\ag-\frac{1}{\dg}D\ag k\cg)(\ag-\frac{1}{\dg}\bg\cg)^{-1}=(D\ag\ag^{-1}\ag-\frac{1}{\dg}D\ag\ag^{-1}\ag k\cg)(\ag-\frac{1}{\dg}\bg\cg)^{-1}$$
which is (as $\ag k=\bg$) $D\ag\ag^{-1}(\ag-\frac{1}{\dg}\bg\cg)(\ag-\frac{1}{\dg}\bg\cg)^{-1}=D\ag\ag^{-1}$. Thus, $D\ag\ag^{-1}$ must be lower triangular.

Conversely if $D\ag\ag^{-1}\in \MLT{n-1}$ and $\bg=\ag k$  for some constant vector $k$ then the  above calculation also shows that $D(\mBB)\mBB^{-1}\in \MLT{n}$. \eproof

\brem{Bstructure} The above result completely describes the structure of the matrices $\mBB$ such that $D\mBB\mBB^{-1}\in\MLT{n}$. If $\mBB=[\ag_{i,j}]$ is a square matrix such that $\bg_j=\ag_j k_j$ for some constant vector $k_j$, where $\ag_j=[\ag_{i,j}]$ is the $j^{th}$ main diagonal block of $\mBB$ and $\bg_j=[\ag_{1,j+1},\ldots,\ag_{j,j+1}]^T$. We assume that $\ag_1\ne0$ then by induction we see that $\mBB$ is such that $D\mBB\mBB^{-1}\in \MLT{n}$ if and only if  $\bg_j=\ag_j k_j$ for some constant vector $k_j$ and $\ag_{j+1,j+1}\ne \cg_j k_j$ where $\cg_j=[\ag_{j+1,1},\ldots,\ag_{j+1,j}]$ for all $j$.
\erem

\brem{mBBdiag} It is easy to see that a similar argument will give another version of \reflemm{DBBlem} to give that
$$D(\mBB)\mBB^{-1} \mbox{ is diagonal}\Leftrightarrow D\ag\ag^{-1} \mbox{ is diagonal and }\bg=\ag k \mbox{ for some constant vector $k$}.$$

Furthermore, in order that $D(\mBB)\mBB^{-1}$ is diagonal we just need  the row vector ( under the main block of $D(\mBB)\mBB^{-1}$, via \mref{MATinverse})  $D\cg(\ag-\frac{1}{\dg}\bg\cg)^{-1}-D\dg\frac{\cg}{\dg}(\ag-\frac{1}{\dg}\bg\cg)^{-1}=0$.  This means $D\cg-\frac{1}{\dg}D\dg\cg=0$ or $D\dg\cg=\dg D\cg$. Equivalently, if $\cg=[\cg_1,\ldots,\cg_j]$ then $\cg_i=c_i\dg$ for $i=1,\ldots,j$ and some constants $c_i$'s.

\erem

We now let $\mA=[a_{i,j}]$ be a matrix such that $\mBB\mA\in \MLT{n}$. We write $\mBB, \mA$ in block form
$$\mBB=\mat{\ag&\bg\\\cg&\dg},\; \mA=\mat{a&b\\c&d}.$$ 

\blemm{BAlem} There is a matrix $\mBB$ such that $\mBB \mA\in \MLT{n}$ and $D\mBB\mBB^{-1}\in \MLT{n}$ if 
\beqno{MATBA}\mbox{there are constant vectors $k_{j}$ such that  $b_j=a_{j,j}k_{j}$  and $(a_j+k_jc_j)^{-1}$ exists for all $j\ge 1$.}\eeq Here, $a_j$'s are the main diagonal blocks of $\mA$,
$k_j=-[k_{1,j},\ldots,k_{j-1,j}]^T$, $b_j=[a_{1,j},\ldots,a_{j-1,j}]^T$   and  $c_j=[a_{j,1},\ldots,a_{j,j-1}]$.

Conversely, if \mref{MATBA} holds and $k_jc_ja_j^{-1}$ is lower triangular (or zero matrix) for all $j\ge1$ then there is $\mBB$ such that $\mBB \mA\in \MLT{n}$ and $D\mBB\mBB^{-1}\in \MLT{n}$.

\elemm

This lemma completely describes the structure of $\mA$ so that the algebraic part (In addition to \mref{Bmatcondthm}) of the conditions  \refcoro{GenMPMatcoro} can be verified. As we will see later, one need some extra algebraic conditions on $a_{i,j}$ for $i\ge j$ to determine  $\mBB$ such that the diagonal entries of $\mA_d=\mL\mBB\mA$ are positive. Also, we need the analytic assumptions on $a_{i,j}$, $i\ge j$ and the Green function $\mG$ determined by $\mA_d$ (and also $\mB$ if it is present and the same boundary conditions) so that \refcoro{GenMPMatcoro} can be used. This reveals an interesting structure of $\mA$: The cross diffusivities of the $i^{th}$ species must be (constant) multiples of the self diffusivity of the $j^{th}$ species if $i<j$ and, meanwhile, if $i>j$ then the $i^{th}$ component will have more freedom locally.

\bproof We argue along the main diagonal blocks. If $\mBB\mA, D\mBB\mBB^{-1}\in \MLT{n}$  then by \reflemm{DBBlem} we must have $\ag b+\bg d=0$. As $\bg= \ag k$ for some constant vector $k$, we see that $b =-dk$. By \mref{MATinverse} $(a_j+k_jc_j)^{-1}$ must exist. This is to say that, by induction, if $\mA=[a_{i,j}]$ be a matrix (of function entries) then we must have that \mref{MATBA} holds.

Conversely, assume that $\mA$ satisfies \mref{MATBA} and $k_jc_ja_j^{-1}$ is lower triangular for all $j\ge1$. We construct $\mBB$ first using only \mref{MATBA}.   For any function $\ag\ne0$, we take $\mBB=[\ag_{i,j}]$ to be a square matrix such that $\bg_j=\ag_j k_j$ where, inductively (starting with $\ag_1=[\ag]$),  $\ag_j=[\ag_{i,j}]$ is the $j^{th}$ main diagonal block of $\mBB$  and $\bg_j=[\ag_{1,j+1},\ldots,\ag_{j,j+1}]^T$ ($j\ge1$). We also take $\cg_j$ to be any row vector and any nonzero function $\dg_j$  such that $(\dg_j-\cg_j k_j)^{-1}$ exists. We see that $D\mBB\mBB^{-1}\in \MLT{n}$. 

However, in order for $\mBB\mA\in \MLT{n}$ we must also have the row vector $c_j=[a_{j,1},\ldots,a_{j,j-1}]$ such that $\ag_ja_j+\bg_jc_j=\ag_j(a_j+k_jc_j)=\llg_j$ is a lower triangular matrix ($a_j$'s are main diagonal blocks of $\mA$). We then define $\ag_j$'s accordingly by choices of $\ag_{i,j}$'s in terms of $a_{i,j}$'s ($\ag_j=\llg_j(a_j+k_jc_j)^{-1}$ for some lower triangular matrix $\llg_j$, such that $\mBB$ is invertible). Thus, $\bg_j=\llg_j(a_j+k_jc_j)^{-1}k_j$. In addition, let $I_j$ denote the $j\times j$ identity matrix,  $\ag_ja_j=\llg_j(a_j+k_jc_j)^{-1}a_j=\llg_j(I_j+k_jc_ja_j^{-1})$ is lower triangular because
$k_jc_ja_j^{-1}$ is. Hence, we can argue inductively (backward) and complete the proof.  \eproof

If $\ag_j$ was predefined then we have to choose $c_j$ such that $a_j+k_jc_j=\ag_j^{-1}\llg_j$ for some lower triangular matrix $\llg_j$. There should be some extra conditions on $k_j$'s and $a_{i,j}$'s ($i>j$) determined by $k_j,a_{i,j}$ for $i\le j$. So, this can be done by induction. Thus, $\ag_{i,j}$ must be such that this system has a solution $c_j$ for given $a_{i,j}, k_i$ and some $\llg_j$).

\brem{BAKrem} We must assume that $\mA,\mB,K$ can be {\em simultaneously}  transformed (or row equivalent) to lower  triangular matrices by the {\em same} matrix $\mBB$ so that there are some connections in their structures. Note that $\mBB$ was determined by $\mA$. Now, suppose that, in block form, $$\mA=\mat{a&dk\\c&d},\;\mB=\mat{\bar{a}&\bar{d}\bar{k}\\ \bar{c}&\bar{d}}$$ for some constant vectors $k,\bar{k}$ then we must have $\ag=\llg(a+kc)^{-1}=\bar{\llg}(\bar{a}+\bar{k}\bar{c})^{-1}$ and $\bg=\ag k=\ag\bar{k}$ so that $\bar{k}=k$ and $(\bar{a}+\bar{k}\bar{c})(a+kc)^{-1}=\LLg$ for some lower triangular matrix $\LLg=\llg^{-1}\bar{\llg}$. That is, $\mA,\mB,K$ must share the same of constant vectors $k_i$'s. Note also that the entries of $\mB_d=\mL\mBB\mB$ and $K_d=\mL\mBB K$ do not have to be positive.

\erem

Again, we introduce the set $\ALT{n}$ of matrix $\mA$ satisfying the structure of \reflemm{BAlem}. That is, the entries $\mA$ verify \mref{MATBA} and $(a_j+k_jc_j)^{-1}$ exist.

In particular, we can take $\ag_1$ to be a nonzero constant and $\cg_j$ to be constant vector. The $\dg_j$'s, $j>1$, can be functions.

Of course, the condition \mref{MATBA} allows more general structure of $\mA$ than \mref{matrowconda} where $\mA$ can be tranformed to lower triangular form by \mref{BBupper} which in fact required that $\frac{k_{i,j}}{k_{i+1,j}}$'s are constants independent of $j>i$ and the matrix $\mBB$ is a constant one. We need only $D\mBB\mBB^{-1}\in \MLT{n}$ here. The matrix $\mBB$ can be a full matrix and the conditions concerning the Green function $\mG$ are obvious if $D(\mL^{-1})=0$.

\subsection{Combining with a change of variables}
We want to generalize \reftheo{GenMPMat} by combining it with a change of variables. Again, we consider the system 
$$\left\{\barr{ll} -\Div(\mA DW) -KW=F &\mbox{in $\Og$,}\\ W=0&\mbox{on $\partial\Og$}\earr
\right.$$
and make a change of variables $W=\ccT v$ for some invertible matrix $\ccT$. We obtain
\beqno{transfsys}\left\{\barr{ll} -\Div(\mA\ccT Dv) -\mA D\ccT Dv-\Div(\mA D\ccT)v -K\ccT v=F &\mbox{in $\Og$,}\\ v=0&\mbox{on $\partial\Og$}\earr
\right.\eeq

We will apply \reftheo{GenMPMat} to this system with $\mA,\mB,K$ are  $\mA\ccT , -\mA D\ccT, (\Div(\mA D\ccT)v +K\ccT)$ respectively.

In order to apply our previous argument to this system and obtain $v>0$ we need that the matrices $\mBB\mA\ccT, -\mBB\mA D\ccT,-\mBB\Div(\mA D\ccT)\in \MLT{n}$ for some matrix $\mBB$ such that $D\mBB \mBB^{-1}\in \MLT{n}$.

We start with the assumption that $\mBB\mA\ccT\in \MLT{n}$ and $\ccT^{-1}D\ccT\in \MLT{n}$. Note that the matrix $K$ must also satisfy that $\mBB K\ccT\in \MLT{n}$.

We see that $\mBB\mA D\ccT=(\mBB\mA\ccT)(\ccT^{-1}D\ccT)\in\MLT{n}$. It is obvious that $\mA\ccT\in\ALT{n}$ implies $D(\mA\ccT)\in\ALT{n}$. So that $D\mA\ccT+\mA D\ccT\in\ALT{n}$. As $\mBB\mA D\ccT\in \MLT{n}$, we can see that $\mBB D\mA\ccT\in\MLT{n}$ (using the same $\mBB$).

Next, $\mBB\Div(\mA D\ccT)=\mBB D\mA D\ccT+\mBB\mA\Delta\ccT$. Because \beqno{DADT}\mBB D\mA  D\ccT=(\mBB \mA\ccT)(\ccT^{-1}\mA^{-1}\mBB^{-1})(\mBB D\mA\ccT)(\ccT^{-1}D\ccT)\eeq and the  matrices on the right are in $\MLT{n}$ ($\ccT^{-1}\mA^{-1}\mBB^{-1}=(\mBB\mA\ccT)^{-1}\in\MLT{n}$). So is the one on the left hand side. We then just need to show that $\mBB\mA\Delta\ccT\in \MLT{n}$.

We have 
$\mBB\mA\Delta\ccT=(\mBB\mA\ccT)(\ccT^{-1}\Delta\ccT)$ so that we need only show that $\ccT^{-1}\Div(D\ccT)\in\MLT{n}$. Since $\ccT^{-1}D\ccT\in\MLT{n}$ we see that $\Div(\ccT^{-1}D\ccT)=D\ccT^{-1}D\ccT+\ccT^{-1}\Delta\ccT\in\MLT{n}$. As $D\ccT^{-1}D\ccT=-(\ccT^{-1}D\ccT)( \ccT^{-1}D\ccT)\in \MLT{n}$ we obtain that $\ccT^{-1}\Delta\ccT\in\MLT{n}$ as desired.

Let $\ccT=\mT^{-1}$. We see that we need to assume that $\mBB\mA\mT^{-1}, D\mT\mT^{-1}\in\MLT{n}$. Let
$$\mA=\mat{a&b\\c&d},\; \mT=\mat{\ag&\bg\\\cg&\dg}.$$

The assumption $ D\mT\mT^{-1}\in\MLT{n}$ implies $\bg=\ag \hat{k}$ for some constant vector $\hat{k}$. Denote $\thg$ by the number $(\dg-\cg\ag^{-1}\bg)^{-1}\ne0$. By \mref{MATinverse} the last column of $\mA\mT^{-1}$ is $$\mat{-a\ag^{-1}\bg\thg+b\thg\\-c\ag^{-1}\bg\thg+d\thg}.$$
Therefore, $\mBB\mA\mT^{-1}\in\MLT{n}$ and $ D\mBB\mBB^{-1}\in\MLT{n}$ if for some constant vector $\bar{k}$, by \reflemm{BAlem}, $-a\ag^{-1}\bg\thg+b\thg=\bar{k}(-c\ag^{-1}\bg\thg+d\thg)$. Using the fact that $\bg=\ag\hat{k}$, this is equivalent to $-a\hat{k}+b=-c\hat{k}\bar{k}+d\bar{k}$ (with a slight abuse of notation here and later in \mref{blockgencond}: $c\hat{k}\bar{k}$ means $\myprod{c,\hat{k}}\bar{k}$). {\em In fact, to make it simple we can choose $\ag,\bg,\cg,\dg$ (i.e. $\mT$) to be constants}.

In order to construct $\mBB$ from $\mA,\mT$, according to \reflemm{BAlem} and \mref{MATinverse} we need to impose further that (as $\bg=\ag \hat{k}$) the following matrix is lower triangular (the easiest case is when it is zero or $c\dg=d\cg$) $$\bar{k}(c-\frac{d}{\dg}\cg)(\ag-\frac{1}{\dg}\bg\cg)^{-1}(\ag-\frac{1}{\dg}\bg\cg)(\ag-\frac{1}{\dg}\bg\cg)^{-1}=\bar{k}(c\dg-d\cg)(\dg\ag-\ag\hat{k}\cg)^{-1}.$$

Thus, we can start with the matrix $\mA$ satisfying this condition and take any matrix $\mT$ such that, in block form, $\bg=\ag \hat{k}$ then there are $\mBB,\mL$  such that $\mBB\mA\mT^{-1},\mL\in \MLT{n}$. If the assumptions of \reftheo{GenMPMat} are verified then we  can conclude that $v=\mT W\gg 0$ if $\mL\mBB F\gg 0$.

We can introduce the term $\mB DW$ in the system. The equation for $v$ will have the extra terms $\mB\ccT Dv, \mB D\ccT v$ and we need that for the same $\mBB$ in the proof the matrices $\mBB\mB\ccT, \mBB\mB D\ccT$ are in $\MLT{n}$. By the same argument, we see that one needs only assume that $\mBB\mB\ccT\in\MLT{n}$. Thus, we can assume that $\mB$ satisfies \mref{blockgencond} (for the same constant vectors $k,\bar{k}$).

Summarizing the above calculations, we now have the following generalization of \reftheo{GenMPMat}
\btheo{GenMPMatT} Assume that $\mA,\mB,\mT$ are square matrices such that theirs main diagonal blocks \beqno{blockgencond}\left\{\barr{l} \mat{a&b\\c&d}\mbox{ satisfy $-a\hat{k}+b=-c\hat{k}\bar{k}+d\bar{k}$, where $\hat{k},\bar{k}$ are some constant vectors,}
\\
\mT=\mat{\ag&\bg\\\cg&\dg} \mbox{ with $\bg=\ag\hat{k}$ for the same vectors $\hat{k}$.}
\earr\right.\eeq
In addition, assume further that $\bar{k}(c\dg-d\cg)(\dg\ag-\ag\hat{k}\cg)^{-1}$
is lower triangular (or zero matrix).

Then there is a matrix $\mBB$ such that $\mBB\mA\mT^{-1},\mBB\mB\mT^{-1}\in \MLT{n}$.

Furthermore,  assume that $\mBB\mA\mT^{-1}=\mL^{-1}\mA_d$, $\mA_d=\mbox{diag}[\hat{a}_{1,1}\cdots,\hat{a}_{n,n}]$ with $\hat{a}_{i,i}>0$ for some $\mL\in \MLT{n}$.

Suppose that  the matrix $\hat{\mC}=\mL\mBB(\Div(\mA D\mT^{-1})+[K-k Id]\mT^{-1})$ be such that $\off{\hat{\mC}}$ is {\em lower triangular} (observe that $\mL\mBB\Div(\mA D\mT^{-1})$ is lower tringular) and its entries are nonegative for some $k>0$ large. Also, suppose that the entries of $\diag{-\hat{\mC}}$ are nonegative.

Denote $\mC=-\mL\mBB\mB\mT^{-1}+\mL\mBB\mA D(\mT^{-1}) + \mL D(\mL^{-1}) \mA_d=[\mathbf{c}_{i,j}]$. Then $\mC$ is lower triangular. Let $\mG=\mbox{diag}[\mG_i]$ be the Green function of the diagonal system (with the homogenous Dirichlet boundary condition) $$-\Div(\mA_dD\zeta)-\diag{\mC}D\zeta+\diag{-\hat{\mC}}.$$ 

Suppose that

\beqno{BmatcondthmT} \myprod{\mathbf{c}_{i,j}(x),D_x\mG_j(x,y)}\ge0\mbox{ a.e $x,y\in\Og$ and $\forall i>j\ge1$}.\eeq

Then the inverse map ${\cal L}^{-1}$ of the map associated to the system 
$$\left\{\barr{ll} -\Div(\mA DW) +\mB DW+k W-KW=F &\mbox{in $\Og$,}\\ W=0&\mbox{on $\partial\Og$.}\earr
\right.$$
is a map from $\{F\in C(\Og,\RR^n)\,:\, \mL\mBB F\gg 0 \}$ into $\{W\in C(\Og,\RR^n)\,:\, \mT W\gg 0 \}$.

\etheo

Compare with \mref{MATBA} in \reflemm{BAlem} when we do not use $\mT$, the condition \mref{blockgencond} here allows more freedom in $a,b$ but some constraints on $c,d$.

Note that $\mL\mBB(\Div(\mA D\mT^{-1}))=\mL\mBB(D\mA D\mT^{-1}+\mA\Div(D\mT^{-1}))$ and $\ccT=\mT^{-1}$ so that $\mL\mBB(\Div(\mA D\mT^{-1})$ is lower triangular by \mref{DADT} and the paragraph follows it.
However, we still have to assume that $\off{\hat{\mC}}$ is lower triangular and the entries of $\off{\hat{\mC}}, \diag{-\hat{\mC}}$ are nonegative  because $K,k Id$ do not satisfy the structure condition \mref{blockgencond} in general  (thee matrices will be on the right hand side of the transformed system). $k>0$ should be large so that ${\cal L}^{-1}$ exists.

We can let $F=\mathbb{G} W$ where $\mathbb{G}$ is a matrix. Obviously, we can write $\mL\mBB F=(\mL\mBB\mG\mT^{-1})\mT W$ to see that  the above theorem implies (setting $v=\mT W$)
\bcoro{GenMPMatcoro} Assume that $\mA,\mB,K$ satisfy the structure described in \reftheo{GenMPMatT}.  Consider the inverse map ${\cal L}^{-1}$ (uhich is supposed to exist) of the map associated to the system 
$$\left\{\barr{ll} -\Div(\mA DW) +\mB DW+k W-KW=F &\mbox{in $\Og$,}\\ W=0&\mbox{on $\partial\Og$.}\earr
\right.$$

Let $\mM\fg=\mathbf{M}\fg$ where $\mathbf{M}$ is a matrix and  such that $\mL\mBB\mathbf{M}\mT^{-1}$ is a positive matrix (i.e. it map $\{v\in\RR^n\,:\, v\gg 0\}$ into itself).
Then,   ${\cal L}^{-1}\mM:\{v\in C(\Og,\RR^n)\,:\, v\gg 0\}\to\{v\in C(\Og,\RR^n)\,:\, v\gg 0\}$.

\ecoro

Again, we recall \refrem{BAKrem} on the connection of the structures of $\mA,\mB,K,\mathbf{M}$. This result should be compared with the example in \refsec{partialcomp} too. {\em It is important to note that the condition that $\mL\mBB\mathbf{M}\mT^{-1}$ is a positive matrix is independent of $\mathbf{M}\mT^{-1}\in \ALT{n}$}.

\brem{blockgencondrem1} The structure of $\mA$ in \mref{blockgencond} generalizes those of matrices in $\ALT{n}$ where $\hat{k}=0$. Again, $k>0$ needs to be sufficiently large such that ${\cal L}^{-1}$ exists. The constant vectors $\hat{k},\bar{k}$ in \mref{blockgencond} must be the same for $\mA,K$ such that the same $\mBB$ can be used for $\mBB\mA,\mBB K\in\MLT{n}$. However, $\mL$ is found based on $\mA$ only.
\erem

\brem{scalarsys} We see that \reftheo{GenMPMat} and its generalizations rely on the usual maximum principles for scalar equations and an use of Green functions. \reflemm{lopez} was not used here and we had to assume that $\mBB$ was a constant matrix. We see that $K$ can be non cooperative. The same effect can happen here if $\mBB$ is not constant and allow $\mG$ in \refcoro{GenMPMatcoro} to be non-cooperative for appropriate $\mL,\mBB,\mT$ via $\mA$. Thus, we can incorporate $k,K$ into $\mG$ in this case.

\erem

\brem{scalarsys-coop} We see that \refcoro{GenMPMatcoro} is reduced to \reflemm{lopez} if $\mA,\mB$ are diagonal matrices. We can take $\mL,\mBB,\mT$ to be the identity matrix (so that the condition concerning the Green functions is automatically satisfied). We then have to assume that $\mG$ is cooperative as in \reflemm{lopez}. Therefore, the introduction of appropriate cross diffusion into the systems will preserve certain maximum principles even when $\mG$ is not cooperative. 

\erem

{\bf On the sign of $D_x\cg(x) D_x\mG(x,y)$:} 

In addition to the above algebraic condition we also have the analytic condition \mref{Bmatcondthm} involving (derivatives of) Green functions. One wishes to find an conditions on $\mA,\mB$ without those using the Green function ($D_x\mG(x,y)$). This turns out to be very hard. At the moment we can say a little bit if we know that the sign of $D_x\hat{\cg}(x)D_x\mG(x,y)$ is fixed for a.e. $x,y\in\Og$

Since the operator is self-adjoint, $\mG$ is symmetric. That is $\mG(x,y)=\mG(y,x)$ and therefore $D_x\mG(x,y)=D_x\mG(y,x)$.

Suppose that $\hat{\cg}=0$ on $\partial\Og$ then (see \cite[Theorem ]{GW})
$$\iidx{\Og}{\myprod{a(x)D_x\mG(x,y),D_x\hat{\cg}(x)}}=\iidx{\Og}{\myprod{a(x)D_x\mG(y,x),D_x\hat{\cg}(x)}}=\hat{\cg}(y).$$

If $\myprod{D_x\hat{\cg}(x),D_x\mG(x,y)}\ge0$ a.e. $x,y\in\Og$  then $\hat{\cg}(y)\ge0$ (resp. $\hat{\cg}(y)\le0$) for a.e. $y\in\Og$.

\subsection{Strongly positivity}

Until now, our theorem \reftheo{GenMPMatT} and its corollaries give us versions of maximum principles showing that if {\em all} components on the right hand side of the system are positive then so are the components of our solutions. These maximum principles are useful for many purposes, however, if we would like to apply the famous Krein-Rutman theorem then there should be another version which is more suitable. 

The generalization is straightforward by combining the arguments leading to \reftheo{GenMPMatT} and \refcoro{GenMPMatcoro}. The condition \mref{BmatcondthmT} is sufficient to guarantee that the right hand side of the system is positive if $\mL\mBB F>0$ ({\em some} of its components are positive) so that \reflemm{lopez} can be applied (if $\mL\mBB[K-k Id]\mT^{-1}$ is {\em cooperative}) to imply $\mT W\gg 0$ ({\em all} of its components are positive). We then have the following theorem whose setting is quite similar to these assertions excepts the definition of $\hat{\mC}$ and the  assumption that the matrix  $\mL\mBB[K-k Id]\mT^{-1}$ is {\em cooperative}. Because we are going to apply \reflemm{lopez} $K$ does not have to satisfy the same structure as that of $\mA,\mB$.

\btheo{GenMPMatTKR} Assume that $\mA,\mB,\mT$ are square matrices satisfy  \mref{blockgencond}.
Then there are matrices $\mBB,\mL$ described in \reftheo{GenMPMatT}

Let
$\hat{\mC}=D(\mL\mBB)\mA D\mT^{-1}+\mL\mBB\mB D\mT^{-1}$. Suppose that  $\hat{\mC}+\mL\mBB(K-k Id)\mT^{-1}$ is {\em cooperative}.

Denote $\mC=-(D(\mL\mBB)\mA\mT^{-1}+\mL\mBB\mB\mT^{-1})=[\mathbf{c}_{i,j}]$.  Suppose that $\mC$ is diagonal (note that $\mC$ is lower triangular).

Consider the inverse map ${\cal L}^{-1}$, which exists if $k$ is sufficiently large, of the map associated to the system 
\beqno{LMsys}\left\{\barr{ll} -\Div(\mA DW) +\mB DW+\kappa(\mL\mBB)^{-1}\mT W+ kW-KW=F &\mbox{in $\Og$,}\\ W=0&\mbox{on $\partial\Og$}.\earr
\right.\eeq Let $\mM\fg=\mathbf{M}\fg$ where $\mathbf{M}$ is a square matrix.

Suppose that  $\mL\mBB\mathbf{M}\mT^{-1}$  is a map from $\{v\in \RR^n\,:\, v>0 \}$ into itself and $\kappa$ is sufficiently large. Then ${\cal L}^{-1}\mM$ maps $\{v\in C(\Og,\RR^n)\,:\, v> 0\}$ into its interior $\{v\in C(\Og,\RR^n)\,:\, v\gg 0\}$.  That is, ${\cal L}^{-1}\mM$ is a {\em strongly positive operator}.

\etheo

As we see in the argument leading to the structure condition \mref{blockgencond}, to make it simple we can choose $\ag,\bg,\cg,\dg$ (i.e. $\mT$) to be constants.

\bproof We now have to rely on  \reflemm{lopez}, which applies to the transformed system whose main parts are diagonal  (see \refrem{lopezrem}),  and cannot merely use the old induction argument and maximum principles for scalar equations as before. 
We will need only that the right hand side is positive (i.e. {\em some} of its components are positive) and obtain the same conclusion. 

Set $F=\mathbf{M}\fg$. As in the proof of \reftheo{GenMPMatT}, let  $v=\ccT^{-1}W$ with $\ccT=\mT^{-1}$.  Arguging as  \refrem{LopezKR}, multiply the system  for $W$  with $\mL\mBB$  to easily see that the equation ${\cal L}W=\mathbf{M}\fg$ can be written as 
\beqno{transfsysa} -\Div(\mA_d Dv) +D(\mL\mBB)\mA D(\mT^{-1}v)+\mL\mBB\mB D(\mT^{-1}v)+\kappa v-\mL\mBB((K-k Id)\mT^{-1}) v=\mL\mBB\mathbf{M}\fg, \eeq
in $\Og$, $v=0$ on $\partial\Og$. We rewrite it as
\beqno{transfsysb} -\Div(\mA_d Dv) -\mC Dv+\kappa v-\hat{\mC} v-\mL\mBB((K-k Id)\mT^{-1})v=\mL\mBB\mathbf{M}\fg, \eeq where $\mC=-(D(\mL\mBB)\mA\mT^{-1}+\mL\mBB\mB\mT^{-1})$, $\hat{\mC}=D(\mL\mBB)\mA D\mT^{-1}+\mL\mBB\mB D\mT^{-1}$ as defined in the theorem.

If we assume that $\mL\mBB\mathbf{M}\fg>0$  then as $\mC$ is diagonal and $\hat{\mC} +\mL\mBB(K-k Id)\mT^{-1}$ is cooperative by our assumptions  then we can apply \reflemm{lopez} (we still need condition \mref{BmatcondthmT}) to have that $v=\mT W\gg 0$ if $\kappa$ is sufficiently large. Because  $\mL\mBB\mathbf{M}\fg>0$ if $\mL\mBB\mathbf{ M}\mT^{-1}\mT \fg>0$, we then  (recalling that $v=\mT W$) assert that ${\cal L}^{-1}\mM$ is a {\em strongly positive operator} (i.e., it maps the cone $\{\fg\in C(\Og,\RR^n)\,:\, \mT\fg>0 \}$ into its interior $\{W\in C(\Og,\RR^n)\,:\, \mT W\gg0 \}$)
if $\mL\mBB\mathbf{M}$  is a map from $\{\fg\in \RR^n\,:\, \mT\fg>0 \}$ into $\{\fg\in \RR^n\,:\, \fg>0 \}$. That is, if $\mL\mBB\mathbf{M}\mT^{-1}$ is strictly positive. The proof is complete. \eproof

Note that $\mC$ is {\em lower triangular} by the calculations after \mref{DADT}. The last assertion of the theorem is the key point which allows us to apply the Krein-Rutman theorem. We see that  $K$ itself {\em may not need be cooperative} in order to imply the assumption that $\mL\mBB[K-k Id]\mT^{-1}$ is {\em cooperative}. Note that we can allow $k=\kappa=0$ and if $\mL\mBB K\mT^{-1}$ is cooperative then we can only assert that the principal eigenvalue of ${\cal L}^{-1}\mM$ is simple with a positive eigenfunction.

\brem{CCrem} Although that we can prove  that $\mC$ is lower triangular and the condition \mref{BmatcondthmT} concerning the Green functions is a complicated one to check. Here, we are no longer need this because we assume that $\mC$ is diagonal. The simplest case is that $\mC=\hat{\mC}=0$ and it is the case if $\mB=0$ and $\mL\mBB,\mT$ are constant matrices (since  $\bar{K}$ can be {\em any} square matrix, we can choose $\bar{K},k$ such that $\mL\mBB(\bar{K}-k Id)\mT^{-1}$ is {\em cooperative}). A bit more general case is that $\hat{\mC}=0$ and $\mC$ is a diagonal matrix. This occurs easily if $\mT$ is a constant matrix (then $\hat{\mC}=0$), and $\mL\mBB$ is (non constant) diagonal (then $\mC$ is diagonal so that $\mathbf{c}_{i,j}=0$ for all $i>j\ge1$). Then, we can take $\mA=(\mL\mBB)^{-1}\mA_d\mT$ and $\mB=(\mL\mBB)^{-1}\mB_d\mT$
with $\mA_d,\mB_d$ are diagonal and the diagonal entries of $\mA_d$ are positive. This provides nontrivial examples for the theorem.
\erem

\brem{hatCrem} $\mT$ can also be a nonconstant matrix too as in \reflemm{LBTspecial}.

\erem

\brem{mTId} We see that $\mathbf{ M}$  can be competitive for $\mL\mBB\mathbf{M}\mT^{-1}$ to be cooperative if we have appropriate $\mL\mBB$ and $\mT$ (in particular, when $\mT=Id$). If $\mT=Id$ then $\hat{\mC}=0$ and we just need that $\mL\mBB(K-k Id)$ is  cooperative and $\mL\mBB\mathbf{M}$ is strictly positive to get that ${\cal L}^{-1}\mM$ is a strongly positive operator. The extreme case is that when $\mL\mBB=Id$ (i.e., $A$ is diagonal) we can combine $\kappa (\mL\mBB)^{-1}=\kappa Id$ into $k Id$ and \reftheo{GenMPMatTKR} is reduced to the classical result \reflemm{lopez} where we need $\kappa+k$ is large.
\erem

{\bf Concerning the construction of $\mBB$.}
We see that, with $\hat{b}=\hat{d}\hat{k}$, we can write $$\mA\mT^{-1}=\mat{\hat{a}&\hat{b}\\\hat{c}&\hat{d}},$$
$$\hat{a}=(a-\frac{1}{\dg}b\cg)(\ag-\frac{1}{\dg}\ag k\cg)^{-1},\; \hat{b}=(\dg-\cg k)^{-1}(-ak+b),\; \hat{c}=(c-\frac{d}{\dg}\cg)(\ag-\frac{1}{\dg}\ag k\cg)^{-1},\; \hat{d}=(\dg-\cg k)^{-1}(d-ck).$$

Thus, we can define $\mBB=\mat{A&B\\C&D}$ according to \reflemm{BAlem} by taking any lower triangular matrix $\llg$ and constant vector $k'$ and let
$B=Ak',\; A=\llg(\hat{a}+\frac{1}{\hat{d}}\hat{b}\hat{c})^{-1}$. 

We can see that once $\mA,\bar{k},k$ is given and satisfies the condition \mref{blockgencond} then we can freely choose $\ag,\cg,\dg$ such that $\bg=\ag k$ to define $\mT$ (see \reflemm{BAlem}). Then we can take any $C,D,k',\llg$ ($k'$ is a constant vector) as long as $B=Ak'$ to define $\mBB$, by \reflemm{BAlem} again. 

{\bf A counter example}: Let $\fg>0$ be the principal eigenfunction of $-\Delta\fg=\llg_1\fg$. For $a,b,k>0$ to be determined we define $W=[\fg,-k\fg]^T$ and $$\mA=\mat{a&0\\0&d},\; \mG=\mat{a\llg_1&0\\-d\llg_1&1}.$$

We see that $W$ is a nonpositive solution of the diagonal system $$-\Div(\mA DW)+kW=\mat{a\llg_1\fg\\-d\llg_1\fg-k\fg}=\mG W.$$ Note that no matter how large $k>0$ is, \reflemm{lopez} can not be applied here because $\mG$ is non cooperative. However, we will show that there is a positive matrix $\ccK$ such that the cross diffusion system $-\Div(\ccK^{-1}\mA\ccK Dw)+kw=\mG w$, enjoys a maximum principle although $\mG$ is non cooperative. This system is equivalent to the diagonal system $$-\Div(\mA Dw)+kw=\ccK\mG\ccK^{-1}w,$$ which satisfies a maximum principle if $\ccK\mG\ccK^{-1}$ is cooperative by \reflemm{lopez}. Thus, the introduction of appropriate cross diffusion recovers the classical maximum principles.

Indeed, the two eigenvectors of $\mG$ are $v_1=[0,1]^T, v_2=[1,\frac{a\llg_1}{1-a\llg_1}]^T$ to the eigenvalues $1,a\llg_1$ respectively. These two vectors are in $\RR_+^2$ if $1>a\llg_1$. We then let $\ccK$ be a positive matrix such that $\ccK v_1=[0,1]^T$ and $\ccK v_2=[1,0]^T$. Then $\ccK\mG\ccK^{-1}$ is obviously cooperative.

\reftheo{GenMPMat} is almost optimal in a sense that if $m=2$ then $\ALT{m}$ describes all matrices $\mA$ for which our method can be used to obtain certain maximum principle. Again, we consider an easy  system
$$\left\{\barr{ll} -\Div(\mA DW) -KW=F &\mbox{in $\Og$,}\\ W=0&\mbox{on $\partial\Og$.}\earr
\right.$$

Assume that $\ccB$ is a function matrix and multiply $\ccB$ to the system to obtain $$\left\{\barr{ll} -\Div(\ccB\mA DW)+D(\ccB)\mA DW -\ccB KW=\ccB F &\mbox{in $\Og$,}\\ W=0&\mbox{on $\partial\Og$.}\earr
\right.$$

In order to use our argument to this system, we suppose that $\ccB\mA\in\ALT{n}$  and prove that $W\gg 0$. 
Then we can apply \reftheo{GenMPMat} with $\mA,\mB, K$ being $\ccB\mA,D(\ccB)\mA, \ccB K$ (resp.).
The matrices $\mA$ and $\mC=-\mL\mBB\mB + \mL D(\mL^{-1}) \mA_d=-\mL\mBB\mB -  D(\mL)\mBB \mA$ are respectively now $$\ccB\mA \mbox{ and
}\mC=-\mL\mBB D(\ccB)\mA - D(\mL) \mBB\ccB \mA=[-D(\mL\mBB\ccB)+\mL D(\mBB)\ccB]\mA=[\mathbf{c}_{i,j}].$$ Of course, we define $\mL,\mA_d$ by $\mBB\ccB\mA=\mL^{-1} \mA_d$.
We look at $\mC=\mL\mBB\ccB D(\mA)-D(\mL\mBB\ccB\mA)+\mL D(\mBB)\ccB\mA$ and we need that it is lower triangular. Since the last two matrices are in $\MLT{n}$ so we just need $\mL\mBB\ccB D(\mA)$ to be so. As $\mL\mBB\ccB D(\mA)=\mBB\ccB\mA\mA^{-1}D(\mA)$ and $\mBB\ccB\mA\in\MLT{n}$ is invertible, if $\mA^{-1}D(\mA)\in \MLT{n}$. 

We must have that the matrices $\ccB\mA, D(\ccB)\mA, \ccB K$ can be simultaneously transformed to lower triangular ones by the same matrix $\mBB$. Note that this is not the same as saying that they can be simultaneously triangularizable.

First, we are going to specify $\mA$ a bit more. If $\mA^{-1}D(\mA)=-D(\mA^{-1})\mA\in \MLT{n}$ then \reflemm{DBBlem}, applied to $\mA^{-1}$, shows that we must have in general that \beqno{genAcond}\mA=\mat{a&-(d-ca^{-1}b)a(a-d^{-1}bc)^{-1}k\\c&d}.\eeq Here, for some matrix $a$, {\em constant} column vector $k$ and number $d\ne0$ $c$ is a row vector. 

Indeed, if in block form (see \mref{MATinverse})$$\mA=\mat{a&b\\c&d} \Rightarrow \mA^{-1}=\mat{(a-d^{-1}bc)^{-1}&-a^{-1}b(d-ca^{-1}b)^{-1}\\
	-d^{-1}c(a-d^{-1}bc)^{-1}&(d-ca^{-1}b)^{-1}}=:\mat{\hat{a}&\hat{b}\\\hat{c}&\hat{d}}.$$
Applying \reflemm{DBBlem} to $\mA^{-1}$ shows that there is a constant vector $\hat{k}=-k$ such that  $$\hat{b}=-\hat{k}\hat{d}\Rightarrow -a^{-1}b(d-ca^{-1}b)^{-1}=(a-d^{-1}bc)^{-1}k = da^{-1}(d Id-a^{-1}bc)^{-1}k.$$

In general ($(d-ca^{-1}b)^{-1}$ is a number), $b=-(d-ca^{-1}b)a(a-d^{-1}bc)^{-1}k$.

Otherwise,  when $n=2$ then $ca^{-1}b=a^{-1}bc$ is a number, we have $(d-ca^{-1}b)^{-1}=(d-a^{-1}bc)^{-1}$. Canceling, we obtain $b=-dk$.
Also, $(a+kc)^{-1}=(a-d^{-1}bc)^{-1}$. Thus, we have that $\mA\in\ALT{n}$.
Hence, if $m=2$ then our result applies if and only if $\mA\in\ALT{2}$.

If $n=2$ then we can take $\ccB=Id$ so that $D(\ccB)\mA=0$ and the result is obvious. 

In general, in block form, let $$\mA=\mat{a&b\\c&d},\; \ccB=\mat{\ag&\bg\\\cg&\dg},\;\mBB=\mat{A&B\\C&D}.$$

As $\mBB\ccB\mA\in \MLT{n}$, we must have $B=Ak$ for some constant vector $k$ and $$A(\ag b+\bg d)+B(\cg b+\dg d)=0\Rightarrow (\ag b+\bg d)+k(\cg b+\dg d)=0.$$
For this, we can choose $\ag+k\cg=\bg+\dg k=0$. Now, in order that $\mBB D(\ccB)\mA\in \MLT{n}$, we must have $A(D\ag b+D\bg d)+B(D\cg b+D\dg d)=0\Leftrightarrow (D\ag b+D\bg d)+k(D\cg b+D\dg d)=0$. The last equation is equivalent to $(\ag Db+\bg Dd)+k(\cg Db+\dg Dd)=(\ag+k\cg)Db+(\bg+\dg k)Dd=0$. This is obvious because $\ag+k\cg=\bg+\dg k=0$.

Thus, \mref{genAcond} completely describes $\mA$ such that \reftheo{GenMPMat} holds when $\mA\not\in\ALT{n}$ but we can find $\ccB$ such that $\ccB\mA\in\ALT{n}$.

Furthermore, we are going to construct $\ccB=\mat{\ag&\bg\\\cg&\dg}$ such that $\ccB\mA\in\ALT{n}$. Again, by \reflemm{BAlem}, in block form, we just need to present a {\em constant} vector $X$ such that $\ag b+d\bg =X(\cg b+\dg d)$ with $b=-(d-ca^{-1}b)a(a-d^{-1}bc)^{-1}k$ for some constant vector $k$ as in \mref{genAcond}.

Obviously, in block form, for given $a,c,d$ and constant $k$ we can always find $\ag,\bg,\dg$ and constant $X$. 
The matrix $\mBB$ is constructed in \reflemm{DBBlem} using $\mA$, and the matrix $\mL$ follows immediately as $\mBB\mA$ is lower triangular.

\newc{\mTT}{\cal T}

\subsection{\reflemm{DBBlem} is necessary}
We conclude this section by presenting an 
example which shows that the situation is very complicated if $\bg_i$'s are functions and that \reflemm{DBBlem} is necessary. 
We will show that for certain choice of function $b$ we cannot have $W>0$ for all solutions of $-\Div(\mA DW)=F$ so that maximum principles similar to our previous ones when $\bg_i$'s are functions cannot be established without the condition $D\mBB\mBB^{-1}\in \MLT{n}$.

Let $\fg$ be the positive principal eigenfunction associated to the eigenvalue $\llg_1$ of the Laplacian. 
Then $W=[\fg,-\fg]^T$ is a nonpositive solution of the full cross diffusion system of two equations $-\Div(\mA DW)=F$ where $\mA,F$ are defined as follows: Let $a,d>0$ be constants and $b,c$ be functions and $$\mA=\mat{a&b\\c&d},\;F=\left[\barr{c} \llg_1(a-b)\fg+DbD\fg\\\llg_1(c-d)-DcD\fg\earr\right],\; \mBB=\mat{1&-\frac{b}{d}\\0&1},\;\mL=\mat{1&0\\-\frac{cd}{ad-bc}&1}.$$

Then  $$\mBB\mA=\mat{a-\frac{b}{d}c&0\\c&d},\;\mBB F=\left[\barr{c} \llg_1(a-b-\frac{b}{d}(c-d))\fg+DbD\fg\\\llg_1(c-d)\fg-DcD\fg\earr\right]=\left[\barr{c} \llg_1(ad-bc)\fg+DbD\fg\\\llg_1(c-d)\fg-DcD\fg\earr\right]$$ and $\mL\mBB\mA$ is diagonal. Futhermore,
\beqno{LBF}\mL\mBB F=\left[\barr{c} \llg_1(ad-bc)\fg+DbD\fg\\\llg_1\left(-cd+(c-d)\right)\fg+[-\frac{cd}{ad-bc}Db-Dc]D\fg\earr\right].\eeq

For simplicity, we choose $c$ such that $\mL$ is a constant matrix (say $\frac{cd}{ad-bc}=k$, a constant) so that the second condition in \mref{bconstmat} is obvious. For any constant $k$, we just let $c=\frac{kad}{kb+d}$ and note that $Dc=\frac{-k^2ad}{(kb+d)^2}Db$. So, $ \left(-\frac{cd}{ad-bc}Db-Dc\right)D\fg=k\left(-1+\frac{kad}{(kb+d)^2}\right)DbD\fg$. We have $$-\Div(\mL\mBB\mA DW)=\mL\mBB F-\mL D\mBB\mA DW$$ since $\mBB$ is no longer a constant matrix.

Next, we want $\mL\mBB F>0$ to verify \mref{bconstmat} (but $W\not>0$ so that \reftheo{GenMPMat} fails). 
If $c(1-d)>d$, by choosing $d>0$ small and $c>0$ large, then  $\llg_1(ad-bc)\fg+DbD\fg,\llg_1\left(-cd+(c-d)\right)\fg>0$. From \mref{LBF}, we see that we just need $$ \left(-\frac{cd}{ad-bc}Db-Dc\right)D\fg =k\left(-1+\frac{kad}{(kb+d)^2}\right)DbD\fg\ge0.$$

We define $b=\fg$ such that $DbD\fg\ge 0$. The above inequality holds if $k\left(-1+\frac{kad}{(kb+d)^2}\right)>0$. This is equivalent to $k(b^2k^2+(2b-a)dk+d^2)<0$. The quadratic $b^2k^2+(2b-a)dk+d^2$ has two positive roots of the same sign. We just need to choose $k>0$ to be between these two roots. With this choice of $b$, we see that $\mBB$ does not satisfies \reflemm{DBBlem}.

\newc{\ccL}{{\cal L}}
\newc{\hccL}{\hat{{\cal L}}}

\newc{\mbg}{\mathbf{m_{bg}}}

\section{Another version of \reftheo{GenMPMatTKR}} \label{newthm}\eqnoset
The definition of $\ccL$ in  \reftheo{GenMPMatTKR} is quite complicated. Besides some clearer conditions and  the differences in the definitions of $\ccL$, $\hccL$, $\mathbf{ M}$ and $\mM$ here we will reformulate it as follows.

\newcommand{\Pcoop}{\mathbf{P_{coop}}}

\btheo{GenMPMatTKRnew} Assume that $\mA,\mB,\mT$ are square matrices and their main blocks matrices  satisfy 
\mref{blockgencond}. In addition, assume further that $\bar{k}(c\dg-d\cg)(\dg\ag-\ag\hat{k}\cg)^{-1}$
is lower triangular (or zero).

Then there is a  matrix $\mBB$ such that $\mBB\mA\mT^{-1},\mBB\mB\mT^{-1}\in \MLT{n}$.
Furthermore,  assume that $\mBB\mA\mT^{-1}=\mL^{-1}\mA_d$, $\mA_d=\mbox{diag}[\hat{a}_{1,1}\cdots,\hat{a}_{n,n}]$ with the functions $\hat{a}_{i,i}>0$ for some $\mL\in \MLT{n}$.

Suppose that  the matrix  $\hat{\mC}=D(\mL\mBB)\mA D\mT^{-1}+\mL\mBB\mB D\mT^{-1}$ is {\em cooperative}.

Denote $\mC=-(D(\mL\mBB)\mA\mT^{-1}+\mL\mBB\mB\mT^{-1})=[\mathbf{c}_{i,j}]$. Suppose that $\mC$ is diagonal.

Consider the inverse map $\hccL^{-1}$, which exists if $k$ is sufficiently large, of the map associated to the system 
\beqno{LMsys1}\left\{\barr{ll} -\Div(\mA DW) +\mB DW+ kW=F &\mbox{in $\Og$,}\\ W=0&\mbox{on $\partial\Og$}.\earr
\right.\eeq   Suppose that there is some positive matrix $\mP_{pos}$ and cooperative matrix $\Pcoop$ such that  for $\kappa$ sufficiently large
\beqno{MPpos} \kappa\mT^{-1}\mP_{pos}^{-1}\mT> \mT^{-1}\mP_{pos}^{-1}(\Pcoop\mT+k(\mL\mBB)).\eeq

Define the matrices $$ 
\mM:=(\mL\mBB)^{-1}[\mP_{pos}+\kappa Id -\Pcoop]\mT-k Id.$$ Then, $\hccL^{-1}\mM$ is a {\em strongly positive operator}.

\etheo

Obviously, from the choices of $\mL,\mBB,\mT$ we can also write \mref{MPpos} as
\beqno{MPposa} \kappa\mT^{-1}\mP_{pos}^{-1}\mT> \mT^{-1}\mP_{pos}^{-1}(\Pcoop\mT+k\mA_d\mT\mA^{-1}).\eeq
On the other hand, $\mM$ depends mildly on $\mA,\mB$ because besides $\mL\mBB,\mT$ we are free to choose $\mA_d,\mB_d$. Also, the conditions of \reftheo{GenMPMatTKRnew} seem to be technical but we will prove in \refcoro{LBTspecCOND1} that they can be all verified by special choices of $\mL\mBB,\mT$ as in \reflemm{LBTspecial} which will be described later.

\bproof
This theorem follows immediately from \reftheo{GenMPMatTKR}. First of all, if $\hat{\mC}$ is cooperative, then for $\hat{\mC}+\mL\mBB(K-k Id)\mT^{-1}$ to be {\em cooperative} (as required in \reftheo{GenMPMatTKR}) we just need that $\mL\mBB(K-k Id)\mT^{-1}$ is cooperative. Thus, we will choose $K=\bar{K}$ such that $\mL\mBB(\bar{K}-kId)\mT^{-1}=\Pcoop$ for some  cooperative matrix  $\Pcoop$. That is, we define
$K=\bar{K}:=(\mL\mBB)^{-1}\Pcoop\mT+ kId$ here.

If we add $\kappa(\mL\mBB)^{-1}\mT W -\bar{K}W$ to both sides of the system \mref{LMsys1} to reduce it to the form \mref{LMsys}. We just need to show that $\hat{W}=\hccL^{-1}\mM W\ge0$ for any $W>0$. Define $\mathbf{M}:=(\mL\mBB)^{-1}\mP_{pos}\mT$. It is clear that from the definition of $\mM$ here we have
$$\hat{W}=\hccL^{-1}\mM W=\ccL^{-1}(\mathbf{M}+\kappa(\mL\mBB)^{-1}\mT-\bar{K})W=\ccL^{-1}\mathbf{ M}W+\ccL^{-1}(\kappa(\mL\mBB)^{-1}\mT-\bar{K})W.$$
The first term on the right hand side is positive if $\ccL^{-1}\mathbf{ M}$ is a positive operator. This is the case because $\mL\mBB\mathbf{ M}\mT^{-1}=\mP_{pos}>0$ and \reftheo{GenMPMatTKR}.

We now see that one needs $\ccL^{-1}(\kappa(\mL\mBB)^{-1}\mT-\bar{K})W=\ccL^{-1}\mathbf{ M}\mathbf{ M}^{-1}(\kappa(\mL\mBB)^{-1}\mT-\bar{K})W>0$. That is, we just need $\mathbf{ M}^{-1}(\kappa(\mL\mBB)^{-1}\mT-\bar{K})>0$. We write it as
$$\mathbf{ M}^{-1}(\kappa(\mL\mBB)^{-1}\mT-\bar{K})=\kappa(\mT^{-1}(\mL\mBB\mathbf{ M}\mT^{-1})\mT)^{-1}- \mathbf{ M}^{-1}((\mL\mBB)^{-1}\Pcoop\mT+kId)>0,$$
which is $\kappa\mT^{-1}\mP_{pos}^{-1}\mT- \mathbf{ M}^{-1}((\mL\mBB)^{-1}\Pcoop\mT+kId)>0$. Hence, the proof is complete as this is exactly \mref{MPpos} which implies that $\kappa\mT^{-1}\mP_{pos}^{-1}\mT> \mathbf{ M}^{-1}((\mL\mBB)^{-1}\Pcoop\mT+kId)$ (by the defintion of $\mathbf{ M}$). \eproof

In order to apply the theorem to prove that ${\cal L}^{-1}\mM$ is strongly positive (so that $\tau>0$), we need to check the key condition  of \reftheo{GenMPMatTKRnew} on $\hat{\mC}$ and assume \mref{MPpos}. That is, 
\beqno{COND1} \hat{\mC}\mbox{ is {\em cooperative} and }\kappa\mT^{-1}\mP_{pos}^{-1}\mT> \mT^{-1}\mP_{pos}^{-1}(\Pcoop\mT+k(\mL\mBB)).
\eeq

First, let us verify the conditions  $\mC,\hat{\mC}$ in the following special cases.

\blemm{LBTspecial} $\mC$ is diagonal in the following cases
\bdes \item[i)] $\mL\mBB$ is a constant matrix and $\mB=0$.
\item[ii)] $\mT$ is a constant matrix and $\mL\mBB$ is diagonal.
\item[iii)] $\mL\mBB, D\mT\mT^{-1}$ are diagonal.
\edes

In addition, $\hat{\mC}=0$ in cases i), ii) and $\hat{\mC}$ is diagonal in case iii).
\elemm

\bproof The proof is easy. Indeed, for i), if $\mL\mBB$ is a constant matrix and $\mB=0$ then $\mC=\hat{\mC}=0$ and \mref{BmatcondthmT} is obvious.
For ii), we have that $\hat{\mC}=0$ and $\mC=-\mB_d$  is diagonal. 

In the same way for iii), as \beqno{Cdiag}\mC=-D(\mL\mBB)(\mL\mBB)^{-1}(\mL\mBB\mA\mT^{-1})D\mT\mT^{-1}-(\mL\mBB\mB\mT^{-1})=-D(\mL\mBB)(\mL\mBB)^{-1}\mA_dD\mT\mT^{-1}-\mB_d,\eeq we see that $\mC$ is diagonal  if $\mL\mBB,D\mT\mT^{-1}$ are.

The last assertion is obvious from the definition of $\hat{\mC}$ as
$$\hat{\mC}=-D(\mL\mBB)(\mL\mBB)^{-1}\mL\mBB\mA\mT^{-1}D\mT\mT^{-1}-\mL\mBB\mB\mT^{-1}D\mT\mT^{-1}.$$
It is easy to see that $D\mT\mT^{-1}=\ag Id$ if and only if $\ag$ is a constant and $\mT=e^{\ag x_1}\ldots e^{\ag x_n}C$ for some constant invertible matrix $C$. \eproof

\brem{hatCrem} In the case iii), we can specify $\mL\mBB,\mB_d$ further such that the diagonal entries of $\hat{\mC}$ are nonnegative as required by \reftheo{GenMPMatTKR}. Indeed,
$$\hat{\mC}=-D(\mL\mBB)(\mL\mBB)^{-1}(\mL\mBB\mA\mT^{-1})D\mT\mT^{-1}-(\mL\mBB\mB\mT^{-1})D\mT\mT^{-1}=-D(\mL\mBB)(\mL\mBB)^{-1}\mA_dD\mT\mT^{-1}-\mB_dD\mT\mT^{-1}$$ so that we can write $\hat{\mC}=C_1+C_2$. Since the (diagonal) entries of $\mA_d$ are positive, the signs of those of $C_1$ are determined by the products of diagonal entries of $-D(\mL\mBB)(\mL\mBB)^{-1}$ and $D\mT\mT^{-1}$. Thus, we can choose $\mL\mBB$ (whose entries are functions) such that these products are nonnegative. Similarly, we choose $\mB_d$ to have that the diagonal entries of $C_2$ are nonnegative.  Hence, we can have that the diagonal entries of $\hat{\mC}$ are nonnegative for any given nonconstant $\mT$.
\erem

In particular, we use consider the cases described in \reflemm{LBTspecial} to have the assumptions on $\mC, \hat{\mC}$ of \reftheo{GenMPMatTKRnew} fulfilled. The only matter left is prove that $\Pcoop, \mP_{pos}$ exist and satisfy \mref{COND1}. We have the following result.

\bcoro{LBTspecCOND1} Assume the cases i)-iii) of \reflemm{LBTspecial}. There are a positive matrix $\mP_{pos}$ and a cooperative matrix $\Pcoop$ such that all conditions of \reftheo{GenMPMatTKRnew} are verified so that $\hccL^{-1}\mM$ is a {\em strongly positive operator}.

\ecoro

\bproof Assume first i). That is $\mL\mBB$ is a constant matrix and $\mB=0$. Let $\mT=cId$ be any diagonal matrix (we can take $\hat{k}=0$ in \mref{blockgencond}, $\mA$ can be still non diagonal). 

In this case $\mT$ commutes with any matrix and \mref{COND1} is reduced to $\mP_{pos}^{-1}[\kappa Id- \Pcoop-c^{-1}k\mL\mBB]>0$. As $\mP_{pos}^{-1}$ is inverse-positive, there is $x>0$ such that $\mP_{pos}^{-1}x>0$. Obviously, for $X=[x|\cdots|x]$ we can find constant matrices $\Pcoop,\mL\mBB$ so that $\kappa Id- \Pcoop-c^{-1}k\mL\mBB=X$ and we can have \mref{COND1}.

For such $\mP_{pos},\Pcoop$ and given $\mbg$ and diagonal$ \mA_d$ (with positive entries) we can also take $\mA=(\mbg+kId)[\mP_{pos}-\Pcoop+\kappa Id]^{-1}\mA_d$.
Hence, every condition of \reftheo{GenMPMatTKRnew} is verified in this case.

On the other hand, suppose ii) or iii) of \reflemm{LBTspecial}. Let $\mL\mBB=c Id$ and $\mT=e^aC$ for some invertible constant matrix $C$ and $e^a=e^{\ag(x_1+\cdots x_n)}$ for some constant $\ag$ and we have case iii) (and $\mA$ can be non-diagonal). If $\ag=0$ then $\mT$ is a conntant matrix and we have case ii).

In this case $\mL\mBB$ commutes with any matrix. By putting $\mL\mBB=c Id, \mT=e^aC$ into \mref{COND1}, we reduce it to $\kappa C^{-1}\mP_{pos}^{-1}C> e^{-a}C^{-1}\mP_{pos}^{-1}(e^a\Pcoop C-kc Id)$. Now, taking $\mP_{pos}$ to be a constant positive matrix and $C=\mP_{pos}^{-1}$, we see easily that \mref{COND1} is equivalent to \beqno{COND1check}(\kappa Id-\Pcoop)\mP_{pos}^{-1}>-e^{-a} kc Id.\eeq

We are going to show that there are $\kappa>0$ and suitable matrices $\Pcoop, \mP_{pos}$ such that \mref{COND1} in the above form holds so that every condition of \reftheo{GenMPMatTKRnew} is verified again. Since the off diagonal entries of the right hand are zero; the diagonal ones are negative if we choose $c>0$, we see that the above holds if we have that the entries of the left hand side matrix are positive.

Applying \reflemm{ZMlem} below to $P=\mP_{pos}^{-1}$, $A=\Pcoop$ and putting $P=sId-B$ into \mref{COND1check} we see that if we choose $\Pcoop, B$ such that their entries are greater that $2$ then $(\kappa Id-\Pcoop)(sId-B)>0$.

Thus, \mref{COND1} is verified by such $\Pcoop,\mP_{pos}$ and
every condition of \reftheo{GenMPMatTKRnew} is verified.  \eproof

In the proof we used the following elementary fact
\blemm{ZMlem}
If $P$ inverse-positive (i.e. $P^{-1}>0$) then $P$ is a $n \times n$ real $Z$-matrix and a non-singular $M$-matrix. This is the case  $P$ is then can be  expressed in the form $P=sId-B$, where $B = [b_{ij}]$ with $b_{ij} \ge 0$, for all $1 \le i,j \le n$, , where $s$ is at least as large as the maximum of the moduli of the eigenvalues of $B$.

In addition, if $t,s>0$, $A=[a_{i,j}],B=[b_{i,j}]$ with $a_{i,j}, b_{i,j}>2$  then $(tId-A)(sId-B)>0$.
\elemm

\bproof The first assertion is a well known fact 
from the theory of positive matrices (e.g, see \cite{BP}). Also, we see that $(t Id-A)(sId-B)=tsId-(A+B)+A B$.
Obviously, if  the entries of $A, B$ are greater than $2$ then $AB>A+B$ and the last assertion follows. \eproof

\brem{MLBTrem} If   $(\mL\mBB)^{-1}(\kappa Id-\Pcoop)\mT-kId>0$ and $\mL\mBB>0$ then these give $\kappa\mT>\Pcoop\mT+k(\mL\mBB)$
so that if  $\mT^{-1}\mP_{pos}^{-1}$ is positive  then from this we imply the most crucial assumption \mref{COND1} of \reftheo{GenMPMatTKRnew} $\kappa\mT^{-1}\mP_{pos}^{-1}\mT> \mT^{-1}\mP_{pos}^{-1}(\Pcoop\mT+k(\mL\mBB))$.

In fact, because $\mP_{pos}^{-1}$ is inverse-positive can write $\mP_{pos}^{-1}=sId-A$ for some positive matrix $A$ (\cite{BP}). If we assume that $\mT^{-1}=tId-B$ then  $\mT^{-1}\mP_{pos}^{-1}$ is positive (for suitable $A,B$, see \reflemm{ZMlem}). Because $\mT^{-1}$ is a $Z$-matrix and non-singular $M$-matrix, this also gives $\mT>0$. For the same reason, we then write $(\mL\mBB)^{-1}=rId-C$ for some positive matrix $C$ so that $$(\mL\mBB)^{-1}\mP_{pos}\mT=(\mL\mBB)^{-1}\mP_{pos}(\mT)=(rId-C)\mP_{pos}\mT=r\mP_{pos}\mT-C\mP_{pos}\mT$$
is positive if $r>0$ large and can be negative if $r,C$ are chosen properly.

Otherwise, if $\mM=(\mL\mBB)^{-1}\mP_{pos}\mT+(\mL\mBB)^{-1}(\kappa Id-\Pcoop)\mT-kId$ is somewhat competitive then we may assume naturally that  $(\mL\mBB)^{-1}(\kappa Id-\Pcoop)\mT-kId<0$. If $\mL\mBB>0$ then this implies $\kappa\mT<\Pcoop\mT+k(\mL\mBB)$
so that if  $-\mT^{-1}\mP_{pos}^{-1}$ is positive  then from this we imply \mref{COND1} $$\kappa\mT^{-1}\mP_{pos}^{-1}\mT> \mT^{-1}\mP_{pos}^{-1}(\Pcoop\mT+k(\mL\mBB)).$$

Again, we can write $\mP_{pos}^{-1}=sId-A$ and $(\mL\mBB)^{-1}=rId-C$ for some positive matrices $A,C$. If we assume that $-\mT^{-1}=tId-B$ then  $-\mT^{-1}\mP_{pos}^{-1}$ is positive (for suitable $A,B$). Because $-\mT^{-1}$ is a $Z$-matrix and non-singular $M$-matrix, this also gives $-\mT>0$. Therefore, $$(\mL\mBB)^{-1}\mP_{pos}\mT=-(\mL\mBB)^{-1}\mP_{pos}(-\mT)=(-rId+C)\mP_{pos}(-\mT)=-r\mP_{pos}(-\mT)+C\mP_{pos}(-\mT)$$
is negative if $r>0$ large so that $\mM$ is completely competitive. Note also that the cases of \reflemm{LBTspecial} can be used here so that all conditions of \reftheo{GenMPMatTKRnew} are satisfied. This should be compared with the examples presented at the end of \refsec{matrix}.
\erem

Finally, the following consequence of \reftheo{GenMPMatTKRnew} which shows that $\hccL$ exists if $\mM$ is properly given.
\bcoro{Mgiven}
Let $\nu_*=\pm1$ and $\mL\mBB$ be satifying one of the cases of \reflemm{LBTspecial} with $\mL\mBB>0$ or its diagonal entries are positive. Suppose that $$\mM=(\mL\mBB)^{-1}\mP_{pos}\mT+(\mL\mBB)^{-1}(\kappa Id-\Pcoop)\mT-kId=(\mL\mBB)^{-1}\mP_{pos}\mT+\nu_*\mM_*,$$ where $\mM_*=\nu_*[(\mL\mBB)^{-1}(\kappa Id-\Pcoop)\mT-kId]$ for some $k,\kappa>0$ and matrices $\mP_{pos}>0$, cooperative $\Pcoop$, and $\mT$  such that $\nu_*\mT>0$ and $D\mT\mT^{-1}\in \MLT{n}$.

If $\mM_*>0$ then there is $\hccL$ as in \reftheo{GenMPMatTKRnew} such that $\hccL^{-1}\mM$ is strongly positive. In addition, we have that $(\mL\mBB)^{-1}\mP_{pos}\mT$
is positive if $\nu_*=1$ and negative if $\nu_*=-1$.

\ecoro

\bproof
In fact, we can write
$\mP_{pos}^{-1}=sId-A$ and $\nu_*\mT^{-1}=t Id- B$ (as $\nu_*\mT>0$ so that $\nu_*\mT^{-1}$  is a real $Z$-matrix and a non-singular $M$-matrix) and assume first that entries of $A,B$ greater than 2 so that $\nu_*\mT^{-1}\mP_{pos}^{-1}$ is positive by \reflemm{ZMlem}. 

Because (recall that $\mL\mBB>0$  or its diagonal entries are positive)  $$\mM_*>0\Rightarrow \nu_*[(\kappa Id-\Pcoop)\mT-k\mL\mBB]>0\Leftrightarrow\nu_*\kappa\mT>\nu_*[\Pcoop\mT+k(\mL\mBB)].$$ As  $\nu_*\mT^{-1}\mP_{pos}^{-1}$ is positive, we  obtain \mref{COND1} from $\nu_*\kappa\mT>\nu_*[\Pcoop\mT+k(\mL\mBB)]$. \reftheo{GenMPMatTKRnew} then applies and provides $\hccL$ such that $\hccL^{-1}\mM$ is strongly positive.

Finally, note that if 
$\mP_{pos}^{-1}=sId-A$ and $\nu_*\mT^{-1}=t Id- B$ then by scaling $\mP_{pos},\mT$  by a scalar $\mu^{-1}$ with $\mu>0$ large then the entries of $A,B$ greater than 2. Accordingly, we scale $\mL\mBB$ by $\mu^{-2}$ and $\kappa,\Pcoop$ by $\mu$.
Thus, if there are parameters and matrices as stated such that
$$\mM=(\mL\mBB)^{-1}\mP_{pos}\mT+(\mL\mBB)^{-1}(\kappa Id-\Pcoop)\mT-kId$$ then the above scalings allows the assumptions we made on $A,B$. Note also that we can write  $(\mL\mBB)^{-1}=rId-C$ for some $r>0$ large (or $\mu$) so that $(\mL\mBB)^{-1}\mP_{pos}\mT=\nu_* r \mP_{pos}(\nu_*\mT)-\nu_*C\mP_{pos}(\nu_*\mT)$
is positive if $\nu_*=1$ and negative if $\nu_*=-1$.
\eproof

\bibliographystyle{plain}

\end{document}